\documentclass[referee, pdflatex,sn-mathphys-num]{sn-jnl}
\usepackage{graphicx}%
\usepackage{multirow}%
\usepackage{amsmath,amssymb,amsfonts}%
\usepackage{amsthm}%
\usepackage{mathrsfs}%
\usepackage[title]{appendix}%
\usepackage{xcolor}%
\usepackage{textcomp}%
\usepackage{manyfoot}%
\usepackage{booktabs}%
\usepackage{algorithm}%
\usepackage{algorithmicx}%
\usepackage{algpseudocode}%
\usepackage{listings}%
\newcommand{\figs}{figs/}  

\usepackage[caption=false,labelformat=simple,listofformat=subsimple]{subfig}
\newcommand{\mylab}[3]{\raisebox{#2}[0mm][0mm]{\makebox[0mm][l]{\hspace*{#1}{#3}}}} % figure superscript 

\theoremstyle{thmstyleone}%
\theoremstyle{thmstyletwo}%

\theoremstyle{thmstylethree}%

\begin{document}

\title[A black-box optimization method with polynomial-based kernels and quadratic-optimization annealing]{A black-box optimization method with polynomial-based kernels and quadratic-optimization annealing}

%%=============================================================%%
%% GivenName	-> \fnm{Joergen W.}
%% Particle	-> \spfx{van der} -> surname prefix
%% FamilyName	-> \sur{Ploeg}
%% Suffix	-> \sfx{IV}
%% \author*[1,2]{\fnm{Joergen W.} \spfx{van der} \sur{Ploeg} 
%%  \sfx{IV}}\email{iauthor@gmail.com}
%%=============================================================%%

\author*[1]{\fnm{Yuki} \sur{Minamoto}}\email{yuki.minamoto@fixstars.com}

\author[2]{\fnm{Yuya} \sur{Sakamoto}}\email{yuya.sakamoto@fixstars.com}

\affil*[1]{\orgname{Fixstars Amplify Corporation}, \orgaddress{\city{Minato-ku}, \postcode{108-0023}, \state{Tokyo}, \country{Japan}}}

\affil[2]{\orgname{Fixstars Corporation}, \orgaddress{\city{Minato-ku}, \postcode{108-0023}, \state{Tokyo}, \country{Japan}}}

%%==================================%%
%% Sample for unstructured abstract %%
%%==================================%%

\abstract{We introduce kernel-QA, a black-box optimization (BBO) method that constructs surrogate models analytically using low-order polynomial kernels within a quadratic unconstrained binary optimization (QUBO) framework, enabling efficient utilization of Ising machines. While the underlying techniques--such as polynomial kernels and surrogate-based optimization--are individually established, their integration in kernel-QA reflects a deliberate design tailored to the challenges of high-dimensional BBO. The proposed method has been evaluated on artificial landscapes, ranging from uni-modal to multi-modal, with input dimensions extending to 80 for real variables and 640 for binary variables. The results demonstrate that kernel-QA is particularly effective for optimizing black-box functions characterized by high-dimensional inputs, showcasing its potential as a robust and scalable BBO approach.}


%%\pacs[JEL Classification]{D8, H51}

%%\pacs[MSC Classification]{35A01, 65L10, 65L12, 65L20, 65L70}

\maketitle

%%%%%%%%%%%%%%%%%%%%%%%%%%%%%%%%%%%%%%%
%%%%%%%%%%%%%%%%%%%%%%%%%%%%%%%%%%%%%%%
%%%%%%%%%%%%%%%%%%%%%%%%%%%%%%%%%%%%%%%
\section{Introduction}\label{intro}

Black-box optimization (BBO) is a promising tool in various sectors where the objectives are expensive to evaluate black-box functions such as complex numerical simulations or experimental measurements. For such optimization problems, a serial optimization method is typically chosen. Here, an evaluation point in the search space is determined at each ``cycle'' (i.e., an iteration of the optimization loop) using information gathered from previous evaluations as shown in Fig.~\ref{fig:typical_flow}. Typically, a model function is constructed in each cycle based on a relatively small dataset consisting of input-output pairs to/from the black-box function (plus some uncertainty information for some cases). Then, an input set that may minimize the model function is obtained as a potentially optimal input during the optimization step (2) in Fig.~\ref{fig:typical_flow}. Subsequently, the black-box function is evaluated with the newly obtained input set, and finally, the new input-output pair is added to the dataset before repeating the next cycle.

\begin{figure}[!t]
  \centerline{
    \subfloat{\includegraphics[trim=0cm 0cm 0cm 0cm, clip=true, width=13.5cm]{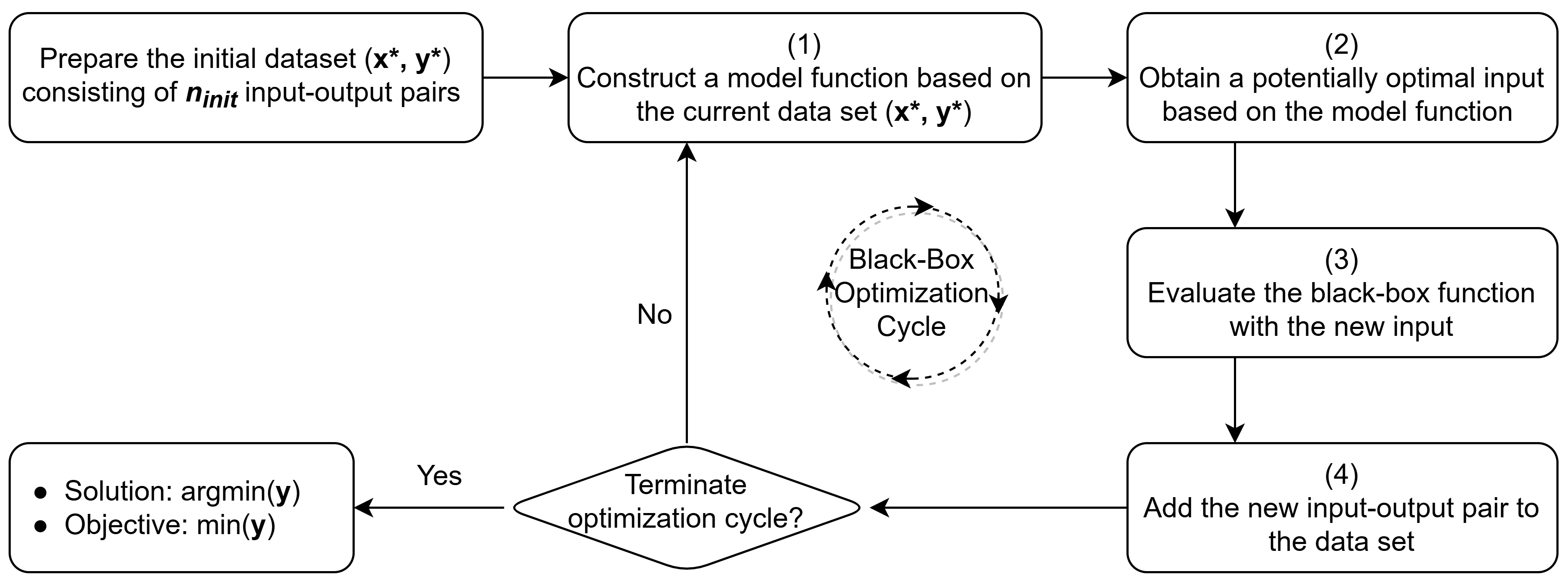}\label{fig:typical_flow}}
  }
  \caption{A flow of a typical serial optimization for a black-box objective function.}%
\label{fig:typical_flow}
\end{figure}

Among existing serial optimization methods, Bayesian optimization stands out for its efficiency in dealing with expensive-to-evaluate black-box functions by building a probabilistic model \citep{frazier2018, Snoek2012, Shahriari2016}. In Bayesian optimization, the Gaussian process (GP) regression typically plays a central role in searching the next evaluation point, balancing the need to explore regions of high uncertainty (to gather more information, typically called ``exploration'') with the need to exploit regions with low expected values (in case of minimization problem as in this study, called ``exploitation''). This exploration makes Bayesian optimization particularly effective for complex black-box functions.

However, the performance of Bayesian optimization degrades as the dimensionality of the search space increases, a challenge known as the curse of dimensionality. This issue arises because the volume of the search space grows exponentially with the number of dimensions, rendering the Gaussian process model less informative. Consequences of high dimensionality include (i) optimization of the model function (step (2) in Fig.~\ref{fig:typical_flow}) becomes increasingly difficult, (ii) high cost in the construction of acquisition function \citep{Baptista2018}, and (iii) insensitivity of the constructed function to the inputs (due to the complexity of the used kernel function) \citep{Hvarfner2024}. As a result, the number of evaluations needed to find an optimal solution becomes prohibitively large, and the method's ability to exploit the model's predictions diminishes.

Factorization machine with quadratic-optimization annealing (FMQA) 
\footnote{FMQA originally stands for factorization machine with {\em quantum} annealing. However, this method is also straightforwardly incorporated with simulated annealers and gate-based quantum computers with a quantum approximate optimization algorithm (QAOA). Thus, a factorization machine with quadratic-optimization annealing may be more appropriate and general, and this point is acknowledged by the authors of the original FMQA study \citep{Kitai2020}.} 
is another serial optimization method for BBO, which has been proposed and applied initially in material informatics \citep{Kitai2020, Inoue2022, Kim2022, Nawa2023}. The key in FMQA is to use factorization machine (FM) \citep{Rendle2010} as a surrogate model for the black-box objective function. The use of an FM model allows the use of Ising machines 
\footnote{An Ising machine is designed to solve combinatorial optimization problems quickly by mimicking the spins' behavior in an Ising model. In an Ising machine, a QUBO problem to be solved is encoded as a network of spins, where each spin represents a binary variable, and the interactions between spins represent the constraints or relationships between these variables. The machine then searches for the spin configuration that minimizes the system's energy, which corresponds to the optimal solution to the problem. Such QUBO solvers include quantum annealers \citep{Dwave} and simulated annealers implemented on GPU or FPGA \citep{Amplify}.}
(simulated or quantum annealing) to solve an optimization problem represented by the model function in a Quadratic Unconstrained Binary Optimization (QUBO) manner at each cycle since a quadratic formulation defines the FM model. FMQA has also been extended and applied to problems with integer/real variables \citep{Izawa2022} and network \citep{Mao2023}. FMQA can handle a relatively large number of decision variables since the model function is optimized using Ising machines. However, the cost of model construction at each cycle is not negligible for large-dimension problems: \citet{Rendle2010} has reported that the cost of training an FM is $\mathcal{O}(n\cdot kd)$, where $n$, $k$ and $d$ are the number of samples in the (current training) dataset, the size of the latent embedding vector (the model hyperparameter), and the number of decision variables (problem dimensions), respectively. Also, stochastic gradient descent typically used to construct an FM may yield locally optimal model coefficients, which is not ideal for a BBO context.

Bayesian optimization of combinatorial structures (BOCS) is a type of Bayesian optimization method where optimization of acquisition function may be performed with simulated or quantum annealing, thereby circumventing the consequence (i) \citep{Baptista2018, Kadowaki2022}. However, the Gibbs sampling used in the model (acquisition function) construction costs $\mathcal{O}(n^2d^2T)$, where $T$ is the number of iterations \citep{Baptista2018}.

This paper addresses high-order black-box optimization (BBO) problems in settings where the number of initial samples is small relative to the input dimensionality. To tackle the curse of dimensionality in such scenarios—particularly in serial optimization—we propose a method called kernel-QA (polynomial-based kernels and quadratic-optimization annealing). While the individual components of this method, such as low-order polynomial kernels and surrogate-based optimization, are well-established, the novelty lies in their deliberate integration: we carefully designed this combination to exploit the strengths of each element in a way that is particularly well-suited for high-dimensional BBO problems. This tailored configuration resolves key challenges (i)–(iii) and leads to a substantial performance gain. We demonstrate the effectiveness of kernel-QA across several benchmark landscapes, where it shows robustness and efficiency, especially in high-dimensional cases where conventional approaches often struggle with computational cost or degraded performance.

The rest of the paper is structured as follows. We describe the formulation of the proposed BBO method in Sec.~\ref{sec:method}. The technique is then assessed using the test functions and conditions described in Secs.~\ref{sec:landscapes} and \ref{sec:conditions}. The results of the assessments are presented and discussed in Sec.~\ref{sec:results}. and important findings in the present study are summarized in the conclusions.

%%%%%%%%%%%%%%%%%%%%%%%%%%%%%%%%%%%%%%%
%%%%%%%%%%%%%%%%%%%%%%%%%%%%%%%%%%%%%%%
%%%%%%%%%%%%%%%%%%%%%%%%%%%%%%%%%%%%%%%
\section{Method}
\label{sec:method}

The keys in kernel-QA are utilizing a QUBO solver (or similar fast optimization solver with/without native constraint support) and the fast construction of the (surrogate) model function for the black-box function. 

%%%%%%%%%%%%%%%%%%%%%%%%%%%%%%%%%%%%%%%
\subsection{Expected value of the black-box function}
\label{sec:expected}

Here, we describe the details of the surrogate model construction for the black-box function (exploitation in a Bayesian optimization context). As noted above, the surrogate model constructed here must be a second-order polynomial to utilize QUBO solvers. Using a relatively low-order polynomial also contributes to circumventing the over-fitting of the constructed model parameters, as the training data is insufficient, which is typical for most BBO problems.

Suppose $f(\mathbf{x})$ is a black-box function with the input variables $\mathbf{x}$. The black-box function $f(\mathbf{x})$ is a complex function (such as numerical simulations or experimental measurements). In BBO, the purpose is to find an input value vector $\mathbf{x}$ that minimizes $f(\mathbf{x})$ with as few cycles as possible (Fig.~\ref{fig:typical_flow}). Let $(\mathbf{x}_i^*, y_i^*)$ is the $i$-th input-output pair to/from $f(\mathbf{x})$, and we have $n$ number of such pairs in the training dataset. Here, $\mathbf{x}_i$ can be a vector of binary, integer, real, and their mix with a dimension $d$. By using appropriate coefficients $\mathbf{Q} \in \mathbb{R}^{d\times d}$ (symmetric matrix), $\mathbf{q} \in \mathbb{R}^{d}$ and a scalar $r$, the output $y_i$ may be written as:
\begin{equation}
    y_i= \mathbf{x}_i^T \mathbf{Q} \mathbf{x}_i + \mathbf{q}^T \mathbf{x}_i + r.
    \label{eq:xQx}
\end{equation}
\noindent
The coefficients, $\mathbf{Q}, \mathbf{q}$, and $r$ can be estimated by solving the following least-squares problem.
\begin{align}
    \hat{\mathbf{Q}}, \hat{\mathbf{q}}, \hat{r} &= \min_{\mathbf{Q}, \mathbf{q}, c} \sum_{i = 1}^{N} \left\{y_i^* - (\mathbf{x}_i^T \mathbf{Q} \mathbf{x}_i + \mathbf{q}^T \mathbf{x}_i + r) \right\}^2 + \lambda \left( \parallel \mathbf{Q} \parallel_F^2 + \parallel \mathbf{q} \parallel_2^2 + r^2 \right).
    \label{eq:least_squares}
\end{align}
\noindent
Here, $\lambda$ is a regularization parameter. Once Eq.~\eqref{eq:least_squares} is solved, we can obtain the surrogate model for the black-box function from Eq.~\eqref{eq:xQx}. However, the number of parameters to optimize for Eq.~\eqref{eq:least_squares} increases with the square of the problem size $d$, which is not ideal for high-dimensional BBO problems.

Now, let our first kernel function be:

\begin{equation}
k_\alpha'(\mathbf{x}_1, \mathbf{x}_2) = (\mathbf{x}_1^{T} \mathbf{x}_2 + \gamma)^2,
\label{eq:kernel_mu}
\end{equation}
\noindent
where $\gamma$ is a constant, and Eq.\eqref{eq:kernel_mu} is essentially a second-order polynomial. Also, let the black-box function $f$ be estimated using the kernel function in Eq.~\eqref{eq:kernel_mu} as:
\begin{align}
    \hat{f}_\mu(\mathbf{x}) &= \sum_{t} c_t k_\mu'(\mathbf{z}_t, \mathbf{x}) \nonumber\\
&= \sum_{t} c_t (\mathbf{z}_t^T \mathbf{x} + \gamma)^2 \nonumber\\
&= \mathbf{x}^{T} (\sum_{t} c_t \mathbf{z}_t  \mathbf{z}_t^T ) \mathbf{x} + 2\gamma (\sum_{t} c_t \mathbf{z})^T \mathbf{x} + c_\mu  \nonumber\\
&= \mathbf{x}^{T} \hat{\mathbf{Q}}_\mu \mathbf{x} + 2\gamma \hat{\mathbf{q}}_\mu^T \mathbf{x} + c_\mu.
\label{eq:bbfunc_kernel}
% \label{eq:mean_linear}
\end{align}
\noindent
Here, all resulting constant terms are put together in $c_\mu$.

We now introduce an optimization problem alternative to Eq.~\eqref{eq:least_squares} to obtain the surrogate model function as:
\begin{align}
    \hat{f}_\mu &= \min_{f} \sum_{i = 1}^N (f(\mathbf{x}_i^*) - y_i^*)^2 + \lambda \parallel f \parallel_{\mathcal{V}}^2 \nonumber \\ 
    &= \min_{\mathbf{c}, \mathbf{Z}} \sum_{i = 1}^N \sum_{t} (c_t k(\mathbf{z}_t, \mathbf{x}_i^*) - y_i^*)^2 + \lambda \sum_{t, t'} c_t c_{t'} k(\mathbf{z}_t, \mathbf{z}_{t'}).
    \label{eq:least_squares2}
\end{align}
\noindent
Here, the subscript $\mathcal{V}$ denotes an operator for the function $f$, $\langle f, f'\rangle_{\mathcal{V}}=\sum_t\sum_t c_t c'_{t'}k_\mu(\mathbf{z}_t, \mathbf{z}'_{t'})$. The optimization problem in Eq.~\eqref{eq:least_squares2} is identical to Eq.~\eqref{eq:least_squares} as $\mathbf{Q}=\sum_t c_t\mathbf{z}_t\mathbf{z}_t^T$, $\mathbf{q} = 2\gamma \sum_{t} c_t \mathbf{z}_t $, $r=c_\mu$. The optimal solution for Eq.~\eqref{eq:least_squares2} is $\hat{f}_\mu(\mathbf{x})=\sum_{j=1}^n c_j k_\mu(\mathbf{x}_j^*,\mathbf{x})$ by the representer theorem without the need for estimating $\mathbf{z}_t$ \citep{scholkopf2001generalized}. Thus, the optimization problem Eq.~\eqref{eq:least_squares2} reduces to the following problem to estimate the optimal $\mathbf{c} \in \mathbb{R}^n$.
\begin{equation}
    \hat{\mathbf{c}} = \min_{\mathbf{c}} \parallel \mathbf{K}^*_\mu\mathbf{c} - \mathbf{y}^*\parallel + \lambda \mathbf{c}^T \mathbf{K}^*_\mu \mathbf{c},
    \label{eq:c_hat_eq}
\end{equation}
\noindent
where $\mathbf{K}^*_\mu$ is a Gram matrix whose $i,j$-element is $k_\mu(\mathbf{x}_i^*, \mathbf{x}_j^*)$. Assuming $\mathbf{K}^*_\mu$ is invertible, Eq.~\eqref{eq:c_hat_eq} can be analytically solved as:
\begin{equation}
    \hat{\mathbf{c}} = (\mathbf{K}^*_\mu + \lambda \mathbf{I})^{-1}\mathbf{y}^*.
    \label{eq:c_hat}
\end{equation}
\noindent
As you can see, the coefficients of the surrogate model are not estimated directly but by kernel regression using a second-order polynomial kernel. When the number of observations is small relative to the dimension of the variables square, it is much cheaper computationally than estimating coefficients of the square order of the dimension. Also, the representer theorem can yield model coefficients that are globally optimal to the samples in the dataset, which is beneficial for black-box optimization.

Note that the inverse operation involved in Eq.~\eqref{eq:c_hat} seems expensive when the number of samples $n$ becomes larger. In the present study, $(\mathbf{K}^*_\mu + \lambda \mathbf{I})^{-1}$ is sequentially updated by using the Woodbury formula \citep{Golub1996} each time a new input-output pair is added at the end of optimization cycle (4) in Fig.~\ref{fig:typical_flow}. This suppresses the computational cost by a factor of $\mathcal{O}(n^2)$, and the overall cost of computing Eq.~\eqref{eq:c_hat} is $\mathcal{O}(n)$.

Finally, the coefficient matrix can be obtained as $\hat{\mathbf{Q}}_\mu = \sum_{i=1}^n \hat{c}_i \mathbf{x}_i^* \mathbf{x}_i^{*T}$, $\hat{\mathbf{q}}_\mu =  \sum_{i=1}^n  \hat{c}_i \mathbf{x}_i^*$, and the final surrogate model $\hat{f}_\mu$ for the black-box function $f$ based on the given dataset is:
\begin{equation}
    \hat{f}_\mu(\mathbf{x}) = \mathbf{x}^{T} \hat{\mathbf{Q}}_\mu \mathbf{x} + 2\gamma\hat{\mathbf{q}}_\mu^T \mathbf{x} + c_\mu.
    \label{eq:mean}
\end{equation}

We can remove the constant $c_\mu$ as this does not generally affect optimization.

%%%%%%%%%%%%%%%%%%%%%%%%%%%%%%%%%%%%%%%
\subsection{Standard deviation}

In some cases, it may be helpful to consider the uncertainty $\hat{f}_\sigma(\mathbf{x})$ of the constructed surrogate model $\hat{f}_\mu(\mathbf{x})$. The general standard deviation of the Gaussian process regression may capture such uncertainty. However, in the context of the kernel-QA method, the formulation of $\hat{f}_\sigma(\mathbf{x})$ must be QUBO-compatible, whereas a well-known acquisition function is not QUBO. Therefore, we construct the acquisition function as QUBO based on the lower confidence bound (LCB):
\begin{equation}
    \hat{f}(\mathbf{x}) = \hat{f}_\mu(\mathbf{x}) - \beta \hat{f}_\sigma(\mathbf{x}),
    \label{eq:acquisition}
\end{equation}
\noindent
where Eq.~\eqref{eq:mean} is used for $\hat{f}_\mu(\mathbf{x})$. Using the Gaussian process regression, the standard deviation $\sigma(\mathbf{x})$ can be written as:
\begin{equation}
    \sigma(\mathbf{x}) = \sqrt{k_\sigma(\mathbf{x}, \mathbf{x}) - \mathbf{k}_\sigma^T(\mathbf{K}^*_\sigma + \lambda \mathbf{I})^{-1} \mathbf{k}_\sigma},
    \label{eq:raw_sigma}
\end{equation}
\noindent
where $\mathbf{K}^*_\sigma$ is another Gram matrix whose $i,j$-element is $k_\sigma(\mathbf{x}^*_i, \mathbf{x}^*_j)$, and $\mathbf{k}_\sigma=[k_\sigma(\mathbf{x}, \mathbf{x}_1^*), k_\sigma(\mathbf{x}, \mathbf{x}_2^*), \cdots, k_\sigma(\mathbf{x}, \mathbf{x}_n^*)]$. Also, for the standard deviation, the following polynomial kernel is used instead of Eq.~\eqref{eq:kernel_mu}.
\begin{equation}
k_{\sigma} (\mathbf{x}_a, \mathbf{x}_b) = \mathbf{x}_a^T \mathbf{x}_b + \gamma.
\label{eq:kernel_sigma}
\end{equation}

The formulation of Eq.~\eqref{eq:raw_sigma} is not a second-order polynomial and therefore requires some modification. Eq.~\eqref{eq:raw_sigma} is simplified to be the QUBO-compatible standard deviation $\hat{f}_\sigma(\mathbf{x})$ as follows.
\begin{align}
    \hat{f}_\sigma(\mathbf{x}) &= k_\sigma(\mathbf{x}, \mathbf{x}) - \mathbf{k}_\sigma^T(\mathbf{K}^*_\sigma + \lambda \mathbf{I})^{-1} \mathbf{k}_\sigma \nonumber\\
&= \mathbf{x}^T \mathbf{x} + \gamma - \sum_{i = 1}^{n}\sum_{j = 1}^{n} L_{ij} ( \mathbf{x}_i^{*T} \mathbf{x} + \gamma)(\mathbf{x}_j^{*T} \mathbf{x} + \gamma) \nonumber\\
&= \mathbf{x}^T \mathbf{x} - \sum_{i = 1}^{n}\sum_{j = 1}^{n} L_{ij} (\mathbf{x}_{i}^{*T} \mathbf{x})(\mathbf{x}_{j}^{*T} \mathbf{x}) - 2 \gamma \sum_{i = 1}^{n}\sum_{j = 1}^{n} L_{ij} (\mathbf{x}_{i}^{*T} \mathbf{x}) + C_v \nonumber\\
&= \mathbf{x}^T \left\{\mathbf{I} - \left(\sum_{i = 1}^{n} \sum_{j = 1}^{n} L_{ij} \mathbf{x}_i^*  \mathbf{x}_j^{*T} \right) \right\} \mathbf{x} - 2 \gamma \sum_{i = 1}^{n}\sum_{j = 1}^{n} L_{ij} (\mathbf{x}_{i}^{*T} \mathbf{x}) + c_\sigma \nonumber\\
&= \mathbf{x}^{T} \hat{\mathbf{Q}}_\sigma \mathbf{x}  - 2 \gamma \hat{\mathbf{q}}_\sigma^T \mathbf{x} + c_\sigma,
    \label{eq:sigma}
\end{align}
\noindent
where $L_{i,j}$ is the $i,j$-element of $(\mathbf{K}^*_\sigma + \lambda \mathbf{I})^{-1}$. In addition, all resulting constant terms are put together in $c_\sigma$. On par with $\hat{f}_\mu$ described in Sec.~\ref{sec:expected}, $(\mathbf{K}^*_\sigma + \lambda \mathbf{I})^{-1}$ is sequentially computed by using the Woodbury formula to suppress computational cost \citep{Golub1996}. Also, as in Eq.~\eqref{eq:mean}, the constant $c_\sigma$ can be removed as this does not affect optimization.

Both $\hat{f}_\mu(\mathbf{x})$ in Eq.~\eqref{eq:mean} and $\hat{f}_\sigma(\mathbf{x})$ in Eq.~\eqref{eq:sigma} are second-order polynomials, thereby the acquisition function in Eq.~\eqref{eq:acquisition} with non-zero $\beta$ being QUBO-compatible. Also, all the coefficients in the polynomial can be analytically obtained rather than iteratively fitted, and this point is advantageous when constructing a surrogate model with a relatively large input dimension.

\subsection{Variable conversion}
\label{sec:encoding}

Kernel-QA inherently requires the variables to be binary to utilise a fast QUBO solver for the model function optimization. However, with appropriate variable encoding, real and integer variables as well as their mix can be considered. Conversion from non-binary to binary variables can be performed either before or after the construction of the model function $\hat{f}(\mathbf{x})$, which are named respectively {\em a priori} and {\em a posteriori} conversions here. 

In the present study, a priori conversion is chosen due to its straightforwardness, with a domain-wall encoding method. A real variable $x$ with lower and upper bounds $(x_{lower}, x_{upper})$ is discretized onto a vector $\mathbf{x}_{disc}$ of a size $n_{bins}$. The spacing of bins $\Delta$ can be either uniform or non-uniform. For example, a variable $x$ with $(x_{lower}, x_{upper})=(0.0, 1.0)$ is discretized into $n_{bins}=5$ bins, and $\mathbf{x}_{disc} = [x_{disc, 1}, \cdots, x_{disc, n_{bins}}] = [0.0, 0.25, 0.50, 0.75, 1.0]$.

Using a binary vector $\mathbf{x}_b$ of a size $n_{bins}-1$, the conversion from $x$ to $\mathbf{x}_b$ is:
\begin{equation}
    x_{b, i} = 
    \begin{cases}
        1, & \text{if }i \le i_{disc}\\
        0, & \text{otherwise}
    \end{cases}
    \:\:(\forall i \in \mathbb{Z} \mid 1 \le i \le n_{bins}-1),
\end{equation}
\noindent
where the discritization index $i_{disc}$ is:
\begin{equation}
i_{disc} = 
    \begin{cases}
        \lfloor(x - x_{lower}) / \Delta + 0.5\rfloor,& \text{if $x$ is uniformly discretized}\\
        \min\{ i \mid x_{disc, i} = x \}, & \text{otherwise}
    \end{cases}.
\end{equation}
\noindent
The conversion from $\mathbf{x}_b$ to $x$ is:
\begin{equation}
    x = \mathbf{x}_{disc}[\sum^{n_{bins}-1}_{i=1} x_{b,i} + 1].
\end{equation}

In the a priori variable conversion method, variable conversion (encoding) from real/integer values to binary vectors is performed before the construction of the surrogate model at (1) in Fig.~\ref{fig:typical_flow}. Thus, any inputs to the constructed surrogate model in Eq.~\eqref{eq:acquisition} are binary values (each non-binary value is converted to a binary vector with the size $n_{bins}-1$). Also, variable conversion (decoding) from binary vectors to real/integer values is required before evaluating the black-box function at (3) in Fig.~\ref{fig:typical_flow}. Evaluation of other encoding methods, such as one-hot encoding and binary encoding with/without {\em a posteriori} conversion, is an interesting topic, which is out of scope in this study.

%%%%%%%%%%%%%%%%%%%%%%%%%%%%%%%%%%%%%%%
\subsection{Exponential transformation of output values}
\label{sec:transform}

In kernel-QA with $\beta=0$, the construction of an appropriate surrogate model (e.g., Eqs.~\eqref{eq:mean}) is the key. Here, the model does not need to be highly accurate: With the complex black-box function and relatively small dataset, such a highly accurate model cannot be obtained, especially for large-dimension problems. In the BBO context, an appropriate model yields a positive correlation with the black-box function. Such correlation characteristics of the model may be assessed by the cross-correlation between the true output values of the black-box function in the training dataset and the values predicted by the trained surrogate model. Typically, such a correlation coefficient should be positive, preferably greater than 0.5 throughout the cycles.

Sometimes, the output values of the black-box function $f(\mathbf{x})$ may have an extensive dynamic range, and the output is sensitive to the input values. Such situations can be identified by looking at the initial data set $(\mathbf{x}^*, \mathbf{y}^*)$. Constructing an ``appropriate'' model for such $f(\mathbf{x})$ can be difficult and often requires some feature scaling on a par with the usual machine learning, for example.

Additional care must be taken for BBO purposes. Minimization problems focus on the model's behavior at relatively small output values. Thus, the model must positively correlate with the black-box function among samples with relatively small output values even though other samples hold large output values. To construct an effective surrogate model for the above-mentioned situation, we consider a transformation that diminishes returns for large values of $f(\mathbf{x})$ and encourage the optimizer to explore regions of the input space where $f(\mathbf{x})$ is smaller. Such transformation is beneficial if there are local minima in the original $f(\mathbf{x})$ that you would like to avoid or de-emphasize in favor of a deeper global search.

One possible transformation would be as follows:
\begin{equation}
    \bar{y} = -\exp{\left(-y/c_m\right)},
    \label{eq:transform}
\end{equation}
\noindent
where $y=f(\mathbf{x})$ and $\bar{y}$ is the transformed output values. Here, $c_m$ is the model parameter, which should be determined based on the initial training dataset. The optimization outcome resulting from this transformation is relatively insensitive to the choice of $c_m$ (see Sec.~\ref{sec:a_exp}). In the present study, we used the ensemble average $\langle \cdot \rangle$ of the output values in the initial training dataset $\mathbf{y}_{init}$:
\begin{equation}
    c_m = \alpha_{exp} \langle \mathbf{y}_{init} \rangle,
    \label{eq:cm}
\end{equation}
where $\alpha_{exp}=1$ by default. We also considered cases with $\alpha_{exp}\ne1$ to assess the sensitivity of the choice of $c_m$ on the overall optimization performance in Sec.~\ref{sec:a_exp}. $c_m$ needs to be a positive number, which may require some linear transformation so that the values $y$ in Eq.~\eqref{eq:transform} are mostly positive before performing Eq.~\eqref{eq:transform}. Perhaps, a straightforward way is to subtract $\min(\mathbf{y}_{init})$ from the original output values $\mathbf{y}_{init}$ if $\min(\mathbf{y}_{init})<0$.

The exponential transformation in Eq.~\eqref{eq:transform} magnifies small output values of $f(\mathbf{x})$ while compressing large ones. Such transformation would be helpful in problems where the objective function has a long tail, or the optimizer needs to focus on further reducing already small values of $f(\mathbf{x})$.

%%%%%%%%%%%%%%%%%%%%%%%%%%%%%%%%%%%%%%%
%%%%%%%%%%%%%%%%%%%%%%%%%%%%%%%%%%%%%%%
%%%%%%%%%%%%%%%%%%%%%%%%%%%%%%%%%%%%%%%

\section{Assessments}
\label{sec:assessments}

Here, we perform various assessments of the proposed kernel-QA by taking the test functions described in Sec.~\ref{sec:landscapes} as black-box with binary and real decision variables.

\subsection{Artificial landscapes}
\label{sec:landscapes}

In the present study, we consider the following artificial landscapes to assess the proposed method, kernel-QA. These functions are for the input of arbitrary dimensions $d$ and are formulated as follows.

\begin{itemize}

\item The Rosenbrock function:
\begin{equation}
f(\mathbf{x}) = \sum_{i=1}^{n-1} \left[ (1 - x_i)^2 + 100 (x_{i+1} - x_i^2)^2 \right],
\label{eq:rosenbrock}
\end{equation}
Global minimum: $f(1, 1, \cdots, 1)=0$.

\item The Rastrigin functions:
\begin{equation}
f(\mathbf{x}) = 10n + \sum_{i=1}^{n} \left[ x_i^2 - 10 \cos(2 \pi x_i) \right],
\label{eq:rastrigin}
\end{equation}
Global minimum: $f(0, 0, \cdots, 0)=0$.

\end{itemize}

\begin{figure}[!t]
  \centerline{
    \subfloat{  \includegraphics[trim=1cm 0cm 6cm 0cm, clip=true, width=6.5cm]{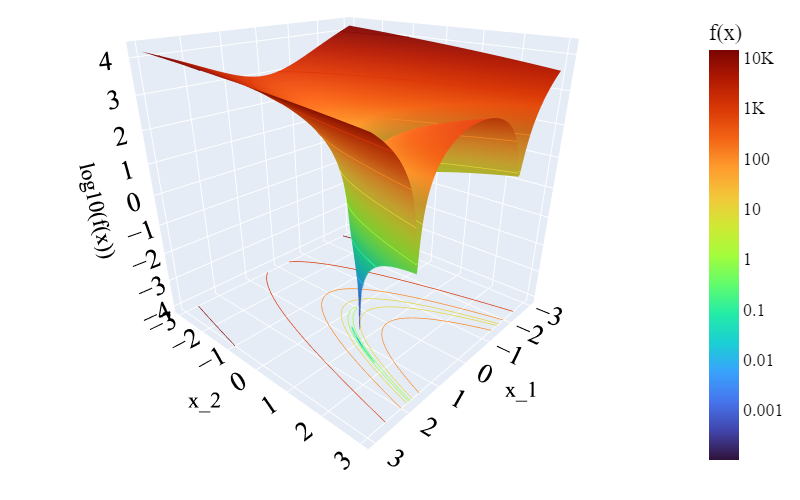}\label{fig:rosenbrock}}
    \mylab{-6.5cm}{5cm}{(a)}
    \subfloat{  \includegraphics[trim=1cm 0cm 6cm 0cm, clip=true, width=7cm]{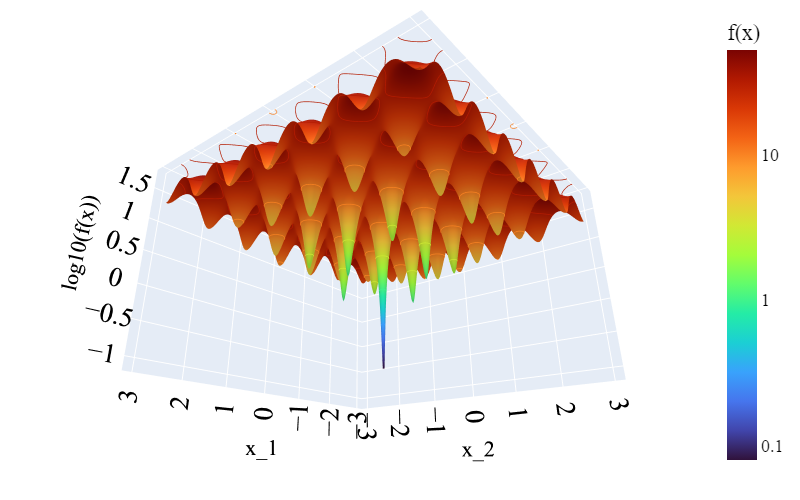}\label{fig:rastrigin}}
    \mylab{-6.5cm}{5cm}{(b)} 
  }
  \caption{The two-dimensional ($d=2$) landscapes of the test functions, (a) Rosenbrock and (b) Rastrigin functions, in logarithm scale.}
\label{fig:landscapes}
\end{figure}

\noindent
Typical two-dimensional ($d=2)$ logarithmic surfaces of these landscapes are shown in Fig.~\ref{fig:landscapes}. These artificial landscapes have different characteristics. Both of the functions are non-convex. The Rosenbrock function is uni-modal, but the global minimum lies in a narrow, parabolic valley \citep{Picheny2013}. As clearly shown in Fig.~\ref{fig:rastrigin}, the Rastrigin function is an example of a non-linear multi-modal function with various local minima.

Depending on the purpose of the assessment, the input dimensions of $d=5$--$80$ are considered for real-variable cases, and $d=40$--$640$ are considered for binary-variable cases. Note that for the kernel-QA, due to the conversion from $d$ real variables (see Sec.~\ref{sec:encoding}), the variables $\mathbf{x}$ considered in the surrogate model (Eq.~\eqref{eq:acquisition}) are binary. The number of such binary variables that the surrogate model and its optimization need to deal with is $d_B=d(n_{bins}-1)$ in the present study, where all the real variables are encoded using the $n_{bins}$ uniform bins. For the binary-variable cases, such variable conversion is unnecessary.

%
%
%

%%%%%%%%%%%%%%%%%%%%%%%%%%%%%%%%%%%%%%%
\subsection{Assessment conditions}
\label{sec:conditions}

The kernel-QA is parametrically assessed in the present study for the artificial landscapes described in Sec.~\ref{sec:landscapes}. For reference, results from typical Bayesian optimization are also compared. Here, optimization conditions considered in the assessments are described.

\begin{table}[!t]
\caption{Summary of assessment conditions for real- and binary-variable problems. $n_{init}=10$, $\alpha_{exp}=1$, and $\beta=0$ unless otherwise noted.}
\label{tab:conditions}
\def\arraystretch{1.5}%  1 is the default, change whatever you need
\begin{tabular}{c|cccc}
\hline
    name & $d$ (real or binary) & $(x_{low}, x_{up})$ & $n_{bins}$ & $d_B$ \\ \hline \hline
    \texttt{r5n} & 5 (real) & $(-3,3)$ & 301 & 1,500 \\ \hline \hline
    \texttt{r5} & 5 (real) & $(-3,3)$ & 61 &  300 \\ \hline
    \texttt{r10} & 10 (real) & $(-3,3)$ & 61 & 600 \\ \hline
    \texttt{r20} & 20 (real) & $(-3,3)$ & 61 & 1,200 \\ \hline
    \texttt{r40} & 40 (real) & $(-3,3)$ & 61 & 2,400 \\ \hline
    \texttt{r80} & 80 (real) & $(-3,3)$ & 61 & 4,800 \\ \hline \hline
    \texttt{b40} & 40 (binary) & - & - & 40 \\ \hline
    \texttt{b80} & 80 (binary) & - & - & 80 \\ \hline
    \texttt{b160} & 160 (binary) & - & - & 160 \\ \hline
    \texttt{b320} & 320 (binary) & - & - & 320 \\ \hline
    \texttt{b640} & 640 (binary) & - & - & 640 \\ \hline
\end{tabular}
\end{table}

A summary of optimization conditions considered in the present assessments is shown in Table~\ref{tab:conditions}. For both kernel-QA and Bayesian optimization, each artificial landscape in Fig.~\ref{fig:landscapes} (with the shown input dimensions, i.e., the number of real/binary variables, $d$) is regarded as black-box and optimized by considering each condition shown in Table~\ref{tab:conditions}. Also, identical initial exploration data are used for kernel-QA and Bayesian optimization cases, whose length (the number of input-output pairs to/from the black-box function) is denoted by $n_{init}$. The initial exploration data are constructed based on randomly chosen input vectors. For each BBO method, optimization runs are repeated ten times independently with a different initial seed for the initial exploration data to collect some statistics. 
%For some conditions with relatively large problem dimensions $d$, initial exploration data of $n_{init}=10$ may be too small. This extreme setting is to demonstrate the effect of $\hat{f}_\sigma(\mathbf{x})$ in Eq.~\ref{eq:acquisition} by considering the non-zero $\beta$ conditions shown in Table~\ref{tab:conditions}. 

For kernel-QA, the same annealing timeout of 5~seconds is considered for all conditions, and Fixstars Amplify Annealing Engine (Amplify AE) \citep{Amplify} is used for the optimization of the acquisition function. Here, Amplify AE is a GPU-based Ising machine that is available to the public. It can accept QUBO problems with up to 256,000 bits. Thus, many typical BBO problems would be solvable regarding decision variable dimensions. Also, for the variable conversion required for the real decision variable cases, the {\em a priori} conversion (Sec.~\ref{sec:encoding}) is applied. These conditions are also applied to the present FMQA runs considered for some assessments. As for the model function for kernel-QA, Eq.~\eqref{eq:acquisition} is considered with $\beta=0$ unless otherwise noted, without including the linear terms ($\gamma=0$ in Eqs.~\eqref{eq:kernel_mu} and \eqref{eq:kernel_sigma}) for simplicity. For all the cases, the regularization parameter $\lambda=1$.

As for Bayesian optimization, the expected improvement (EI) is used for the acquisition function whose coefficients are determined analytically \citep{Jones1998}. Also, the minimization of the acquisition function is performed by using the limited-memory Broyden-Fletcher-Goldfarb-Shanno (L-BFGS) algorithm \citep{Liu1989}, which would cost $\mathcal{O}(dn^2T)$, where $T$ is the number of iterations required. These settings are typical in Bayesian optimization and implemented using the GPyOpt library \citep{GPyOpt} in the present study. Note that of all the parameters shown in Table~\ref{tab:conditions}, only $d$, $(x_{low}, x_{up})$ and $n_{init}$ are directly relevant to Bayesian optimization.

Ten independent runs (different initial data samples) are performed for each optimization case, and their average and standard deviation are discussed in the present assessment. Note that identical initial training data samples are used between different optimization methods.

%%%%%%%%%%%%%%%%%%%%%%%%%%%%%%%%%%%%%%%
%%%%%%%%%%%%%%%%%%%%%%%%%%%%%%%%%%%%%%%
%%%%%%%%%%%%%%%%%%%%%%%%%%%%%%%%%%%%%%%
\section{Results}
\label{sec:results}

\begin{figure}[!t]
\hspace{-0.5cm}
  \centerline{
    \subfloat{  \includegraphics[trim=0cm 0cm 0cm 0cm, clip=true, height=4.5cm]{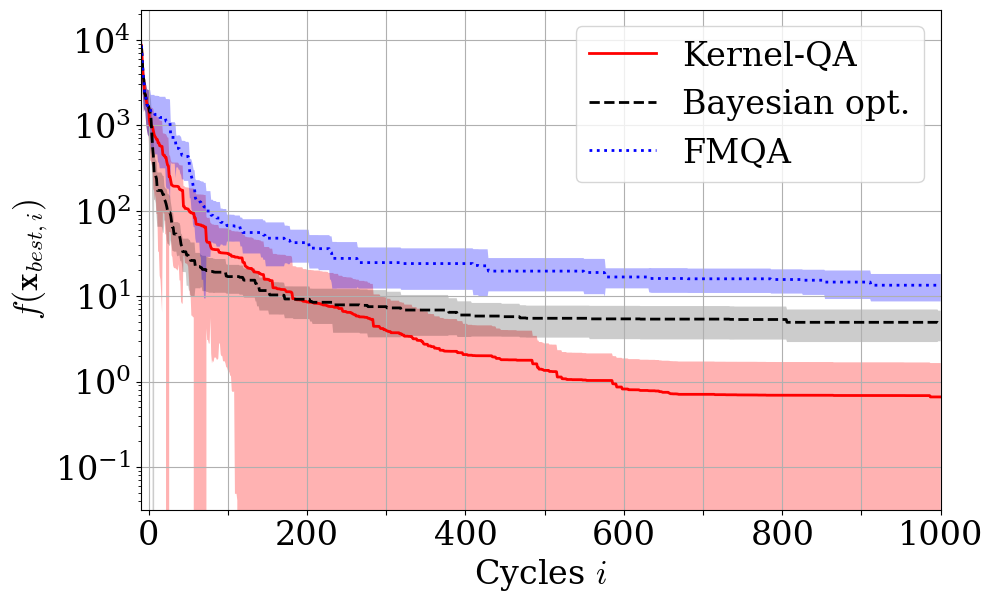}\label{fig:rosenbrock_kq_bo_fm_d5}}
    \mylab{-3.5cm}{-0.5cm}{(a)}
    \hspace{-0.4cm}
    \subfloat{  \includegraphics[trim=0cm 0cm 0cm 0cm, clip=true, height=4.5cm]{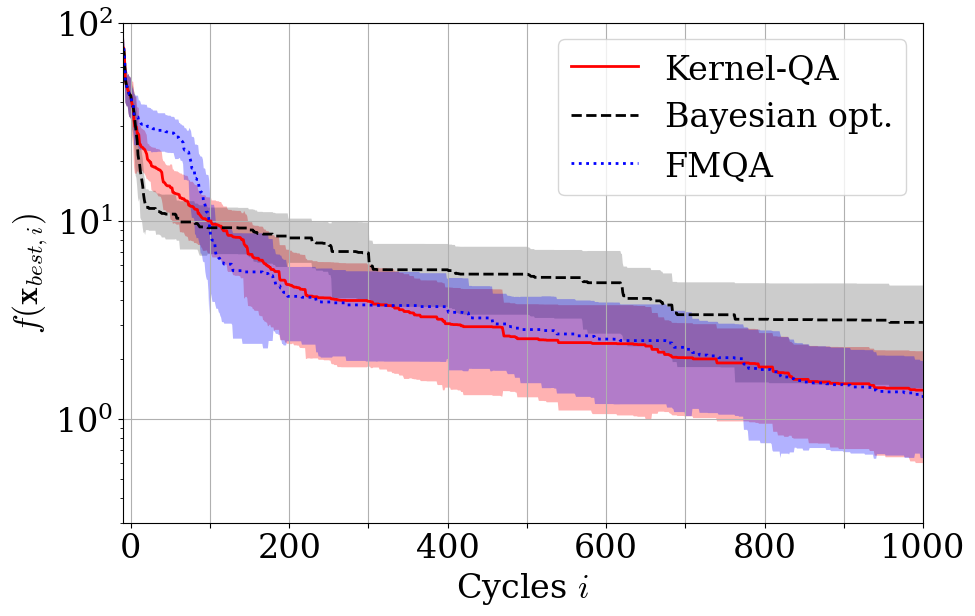}\label{fig:rastrigin_kq_bo_fm_d5}}
    \mylab{-3.5cm}{-0.5cm}{(b)}
  }
  \vspace{0.6cm}
  \caption{Evolution of $f(\mathbf{x}_{best, i})$ for (a) Rosenbrock and (b) Rastrigin functions with the real-variable dimension $d=5$ (\texttt{r5h}). Optimization uses kernel-QA (red) and Bayesian optimization (black). Corresponding results of FMQA are also shown in blue as a reference. Note that plots in the negative cycles evaluate $f(\mathbf{x})$ for the initial training dataset.}
\label{fig:d5}
\end{figure}
\begin{figure}[!t]
  \centerline{
    \subfloat{  \includegraphics[trim=0cm 0cm 0cm 0cm, clip=true, height=4.5cm]{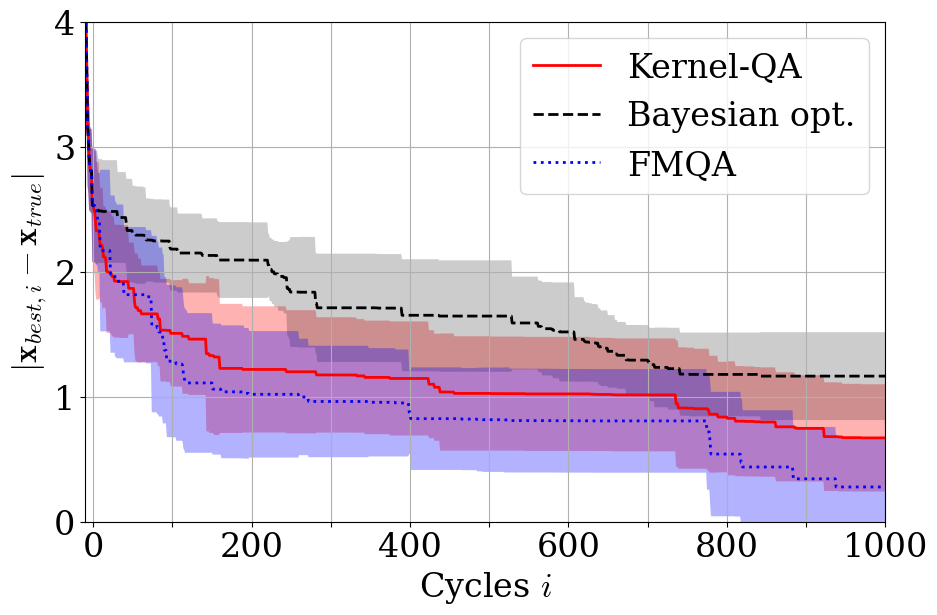}}
  }
  \caption{Evolution of the distance between the found and true solutions obtained for Fig.~\ref{fig:rastrigin_kq_bo_fm_d5}.}
\label{fig:d5_distance}
\end{figure}
%
%
%

%%%%%%%%%%%%%%%%%%%%%%%%%%%%%%%%%%%%%%%
\subsection{General features}
\label{sec:general}

Kernel-QA and Bayesian optimization are performed to find an optimal input for ``black-box'' functions, the Rosenbrock and Rastrigin functions, with the real-number input dimension $d=5$ under the condition labeled as \texttt{r5h} in Table~\ref{tab:conditions}. The evolution of the objective function values with the found best solution at $i$-th cycle, $f(\mathbf{x}_{best, i})$, averaged over ten independent optimization runs, is shown in Fig.~\ref{fig:d5}. Generally, all optimization methods show reasonable trends. That said, the results clarify the different characteristics of different optimization methods.

For the Rosenbrock function (no local minima), after the 1000-th cycle, the values of $f(\mathbf{x}_{best, 1000})$ averaged over the ten runs are 0.7 for kernel-QA, 4.9 for Bayesian optimization, and 13.4 for FMQA. Clearly, FMQA shows a digit worse result, and this is potentially because the stochastic gradient descent used in FM with the relatively small dataset might have resulted in local optimal model coefficients as discussed in Eq.~\eqref{eq:c_hat}, and such behavior does not seem beneficial for the Rosenbrock's landscape, where the global minimum lies in a very narrow valley.

As for the Rastrigin function in Fig.~\ref{fig:rastrigin_kq_bo_fm_d5}, Bayesian optimization yields a rapid decrease of $f(\mathbf{x}_{best, i})$ at the early phase of optimization cycles (say first 100 cycles), but the slope of evolution becomes less than the other methods. Kernel-QA and FMQA catch up with Bayesian optimization after approximately 100 cycles and continue to find better solutions that yield smaller $f(\mathbf{x})$. The values of $f(\mathbf{x}_{best, 1000})$ averaged over the ten runs are 1.4 for kernel-QA, 3.1 for Bayesian optimization, and 1.3 for FMQA.

An interesting observation is shown in Fig.~\ref{fig:d5_distance} which plots the distance between the found best solution at $i$-th cycle $\mathbf{x}_{i, best}$ and the true solution $\mathbf{x}_{true}$, computed as $\|\mathbf{x}_{i, best}-\mathbf{x}_{true}\|$ for the Rastrigin case. The evolution shows that the solutions found by kernel-QA are actually closer to the true solution than the ones found by Bayesian optimization as soon as the optimization cycles are started. On the other hand, in terms of the $f(\mathbf{x}_{best, i})$ shown in Fig.~\ref{fig:rastrigin_kq_bo_fm_d5}, Bayesian optimization shows a rapid decrease and seems outperforms kernel-QA for the first 100 cycles. This conflicting observation reveals kernel-QA and FMQA successfully avoid being trapped in local minima for the Rastrigin function. Such a tendency to avoid local minima in kernel-QA is considered due to the low-order polynomial in its surrogate model, which would circumvent the overfitting of the model parameters with a relatively small dataset typical in BBO problems. A similar trend is shown in Figs.~\ref{fig:rastrigin_kq_bo_fm_d5} and \ref{fig:d5_distance} for FMQA, which also utilizes a second-order polynomial in the surrogate model, ensures this insight.

%%%%%%%%%%%%%%%%%%%%%%%%%%%%%%%%%%%%%%%
\subsection{Mitigation of curse of dimensionality}
\label{sec:curse}

\begin{figure}[!t]
\hspace{-0.4cm}
  \centerline{
    \subfloat{  \includegraphics[trim=0cm 0cm 0cm 0cm, clip=true, height=4.4cm]{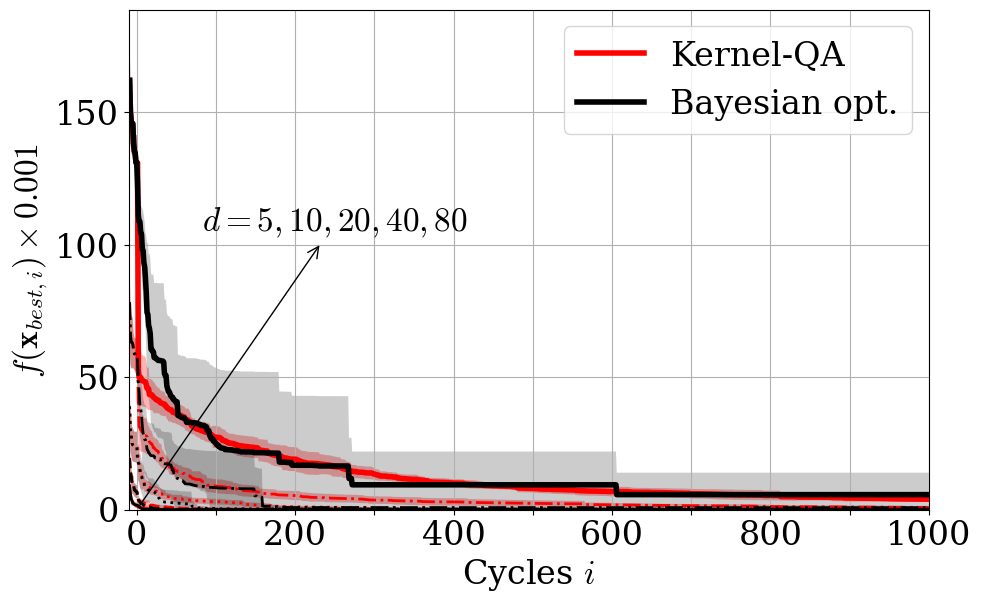}\label{fig:rosenbrock_kq_bo_d5_d80}}
    \mylab{-3.5cm}{-0.5cm}{(a)}
    \hspace{-0.4cm} 
    \subfloat{  \includegraphics[trim=0cm 0cm 0cm 0cm, clip=true, height=4.4cm]{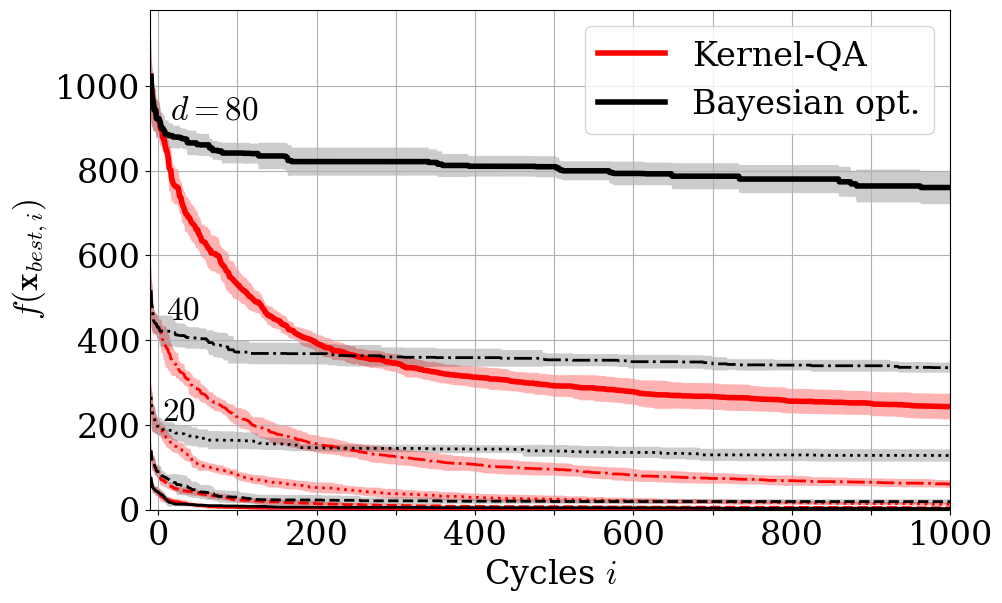}\label{fig:rastrigin_kq_bo_d5_d80}}
    \mylab{-3.5cm}{-0.5cm}{(b)} 
  \vspace{1cm}
  }
  \vspace{0.6cm}
  \caption{Evolution of $f(\mathbf{x}_{best, i})$ with various real-variable dimensions $d=5$ (thin solid), $10$ (dashed), $20$ (dotted), $40$ (dash-dotted) and $80$ (solid line), for (a) Rosenbrock and (b) Rastrigin functions. Optimization uses kernel-QA (red) and Bayesian optimization (black). Note that plots in the negative cycles evaluate $f(\mathbf{x})$ for the initial training dataset.}
\label{fig:d5_d80}
\end{figure}

\subsubsection{Real-variable problems}
\label{sec:real}

Figure~\ref{fig:d5_d80} shows the evolution of $f(\mathbf{x}_{best, i})$ with different real-variable dimensions $d=5$ (\texttt{r5}), 10 (\texttt{r10}), 20 (\texttt{r20}), 40 (\texttt{r40}) and 80 (\texttt{r80}) to see the effect of ``curse of dimensionality'' for kernel-QA and Bayesian optimization.

For the ``black-box'' function with the uni-modal feature (Rosenbrock), kernel-QA and Bayesian optimization perform reasonably at a similar level in an average sense. However, the standard deviation observed for Bayesian optimization is substantially more significant than for kernel-QA throughout the cycles (see Fig.~\ref{fig:rosenbrock_kq_bo_d5_d80}. The large standard deviations for Bayesian optimization become obvious when $d\ge20$ in the present study. For instance, these deviations are 3.2 ($d=20$), 4.4 ($d=40$), and 7.3 ($d=80$) times greater than kernel-QA at the $50$-th optimization cycle. This observation suggests that the quality of the found best solution, especially for cycles $10<i<300$, varies significantly from run to run for the present Bayesian optimization. 

For multi-modal function (Rastrigin), the difference between the two optimization methods is also apparent: optimization progresses at a greater rate and a more consistent manner for kernel-QA as shown by much lower  $f(\mathbf{x}_{best, i})$ values throughout the cycles in Fig.~\ref{fig:rastrigin_kq_bo_d5_d80}. On the other hand, the results show that the adverse effect of increased dimensionality on typical Bayesian optimization is emphasized when the objective function yields multi-modal.

For $d=[5, 10, 20, 40, 80]$, the minimum objective function values found at the end of shown cycles, averaged over the ten optimization runs, are $[1.1, 4.8, 89.4, 893.1, 3950.0]$ (Rosenbrock) and $[1.6, 4.5, 13.0, 60.5, 243.3]$ (Rastrigin) for kernel-QA. As for Bayesian optimization, they are $[5.3, 53.7, 278.7, 793.9, 5700.5]$ (Rosenbrock) and $[2.7, 19.2, 127.9, 335.2, 759.9]$ (Rastrigin).

When considering real-variable problems, Bayesian optimization seems reasonable for uni-modal landscapes like the Rosenbrock function with input dimensions $d\lessapprox20$. At the same time, kernel-QA performs well for optimization problems with the input dimensions $d \gg 20$ for uni- and multi-modal landscapes.

\begin{figure}[!t]
\hspace{-0.4cm}
  \centerline{
    \subfloat{  \includegraphics[trim=0cm 0cm 0cm 0cm, clip=true, height=4.4cm]{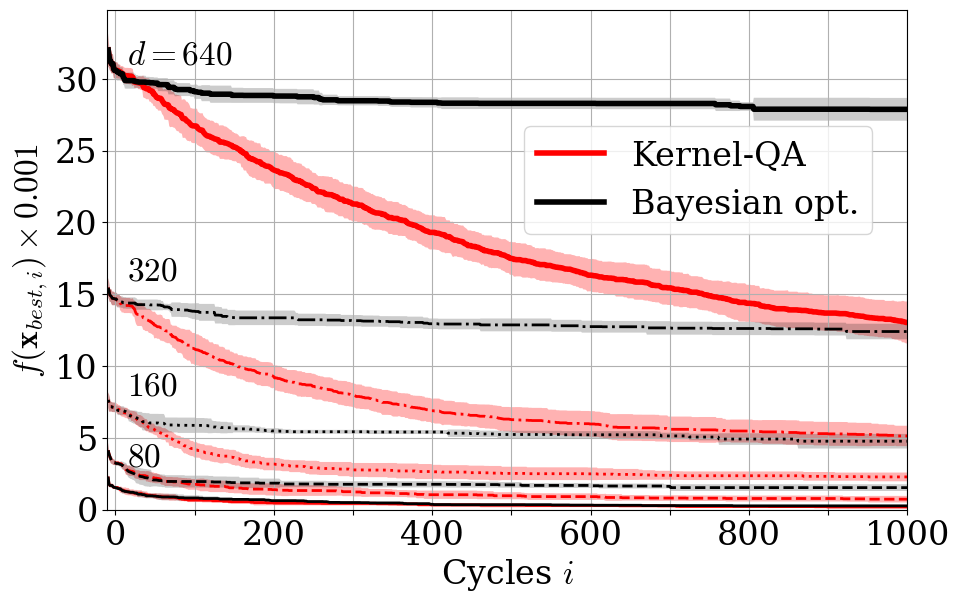}\label{fig:rosenbrock_binary_kq_bo_d40_d640}}
    \mylab{-3.5cm}{-0.5cm}{(a)}
    \hspace{-0.4cm} 
    \subfloat{  \includegraphics[trim=0cm 0cm 0cm 0cm, clip=true, height=4.4cm]{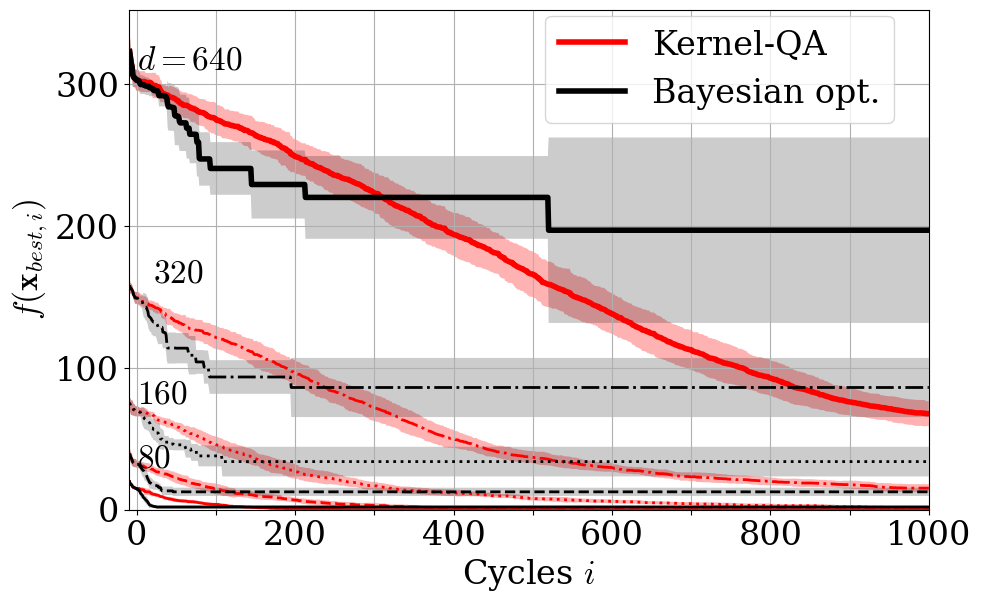}\label{fig:rastrigin_binary_kq_bo_d40_d640}}
    \mylab{-3.5cm}{-0.5cm}{(b)} 
  \vspace{1cm}
  }
  \vspace{0.6cm}
  \caption{Evolution of $f(\mathbf{x}_{best, i})$ with different binary-variable dimensions $d=40$ (thin solid), $80$ (dashed), $160$ (dotted), $320$ (dash-dotted), $640$ (dash-dotted), for (a) Rosenbrock and (b) Rastrigin functions. Optimization uses kernel-QA (red) and Bayesian optimization (black). Note that plots in the negative cycles evaluate $f(\mathbf{x})$ for the initial training dataset.}
\label{fig:d40_d640}
\end{figure}

\subsubsection{Binary-variable problems}
\label{sec:binary}

Given that kernel-QA utilizes the Ising machine, the advantage of kernel-QA for large-dimension problems is even more evident when the optimization cases with binary decision variables (combinatorial BBO) are considered. The test functions are based on Eqs.~\eqref{eq:rosenbrock} and \eqref{eq:rastrigin}. However, half of the input elements $x_i$ are randomly flipped to be $\hat{x}_i$, and $\hat{\mathbf{x}}$ is used as the input to the function.

\begin{equation}
    \hat{x}_{i} = 
\begin{cases}
1 - x_i, & \text{if } i \in \{j_1, j_2, \dots, j_{d/2}\}, \\
x_i, & \text{otherwise},
\end{cases}
\end{equation}
\noindent
where $j_1, j_2, \dots, j_{d/2}$ are randomly chosen unique indices between 1 and $d$. This additional flipping treatment is essential for a fair and meaningful comparison: The used Ising machine, Amplify AE, searches the solution around $\mathbf{x}=\mathbf{0}$ initially, which results in the superior performance of kernel-QA for the problems where $\mathbf{x}_{true} = \mathbf{0}$, but such behavior does not reflect general performance features of kernel-QA.

Figure~\ref{fig:d40_d640} shows assessment results with the binary input dimensions $d=40$ (\texttt{b40}), 80 (\texttt{b80}), 160 (\texttt{b160}), 320 (\texttt{b320}), and 640 (\texttt{b640}). For these input dimensions, the minimum objective function values found at the end of shown cycles (averaged over the ten independent optimization runs) are $[186.1, 734.5, 2284.0, 5136.7, 13039.2]$ (Rosenbrock), and $[0.0, 0.0, 1.2, 15.2, 67.7]$ (Rastrigin) for kernel-QA. As for Bayesian optimization, they are $[270.3, 1538.6, 4769.8, 12407.1, 27876.9]$ (Rosenbrock), and $[1.9, 12.6, 33.9, 86.1, 196.9]$ (Rastrigin). While Bayesian optimization is not particularly known for being advantageous for discrete or binary variables, kernel-QA consistently demonstrates robust and reliable performance even for huge variable dimensions. Also, on par with the real-variable problems discussed in Sec.~\ref{sec:real}, relatively large standard deviations of $f(\mathbf{x}_{best, i})$ are also observed for Bayesian optimization at larger $d$ conditions as shown in Fig.~\ref{fig:rastrigin_binary_kq_bo_d40_d640}.

%%%%%%%%%%%%%%%%%%%%%%%%%%%%%%%%%%%%%%%
\subsection{Effect of $\alpha_{exp}$}
\label{sec:a_exp}

\begin{figure}[!t]
  \vspace{0.5cm}
  \centerline{
    \hspace{-0.7cm}
    \subfloat{  \includegraphics[trim=0cm 0cm 0cm 0cm, clip=true, height=4.4cm]{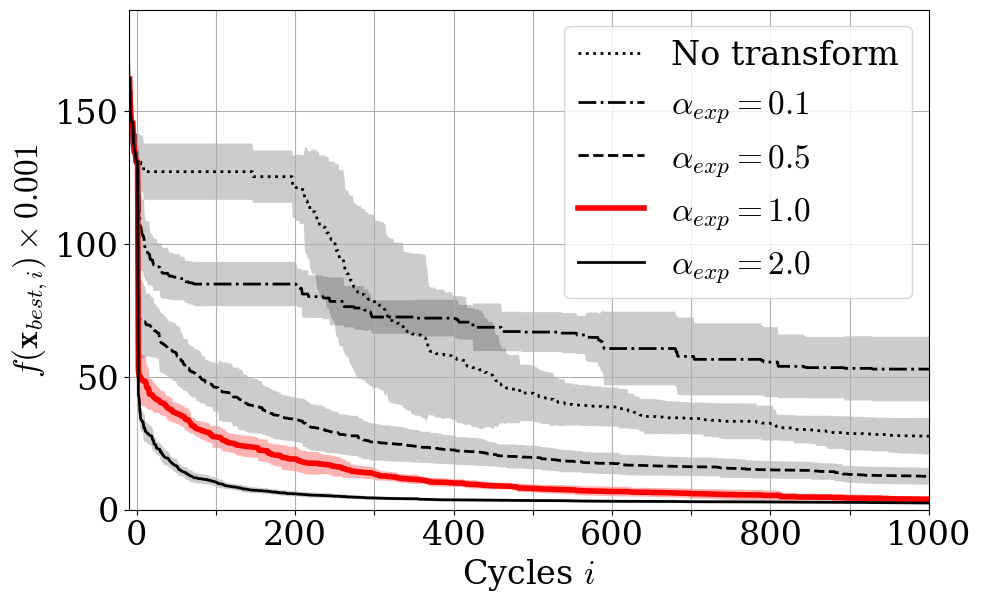}\label{fig:rosenbrock_kq_d80_exp}}
    \mylab{-3.5cm}{-0.5cm}{(a)}
    \hspace{-0.4cm} 
    \subfloat{  \includegraphics[trim=0cm 0cm 0cm 0cm, clip=true, height=4.4cm]{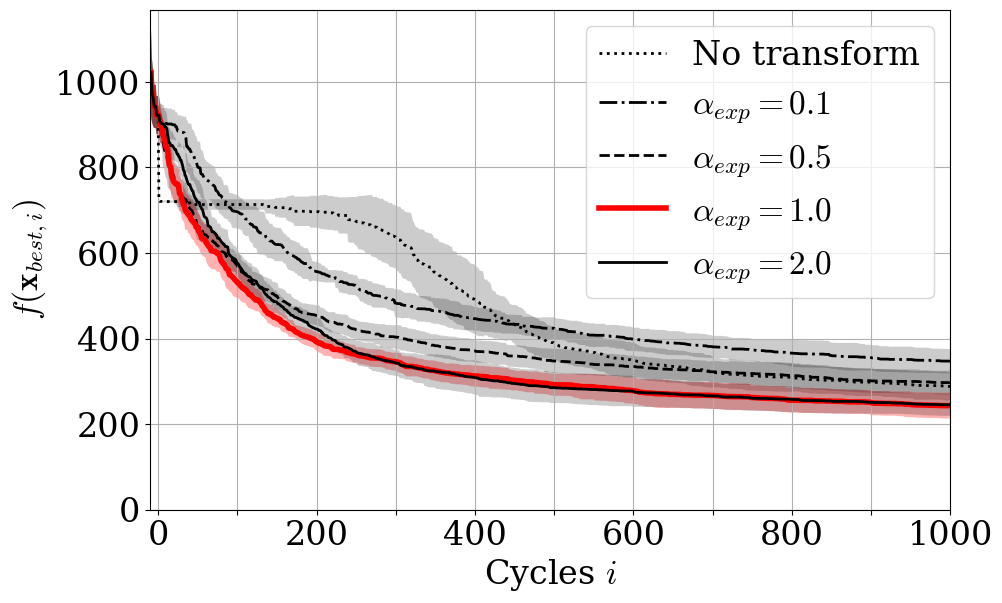}\label{fig:rastrigin_kq_d80_exp}}
    \mylab{-3.5cm}{-0.5cm}{(b)} 
  }
  \vspace{0.5cm}
  \centerline{
    \hspace{-0.6cm}
    \subfloat{  \includegraphics[trim=0cm 0cm 0cm 0cm, clip=true, height=4.4cm]{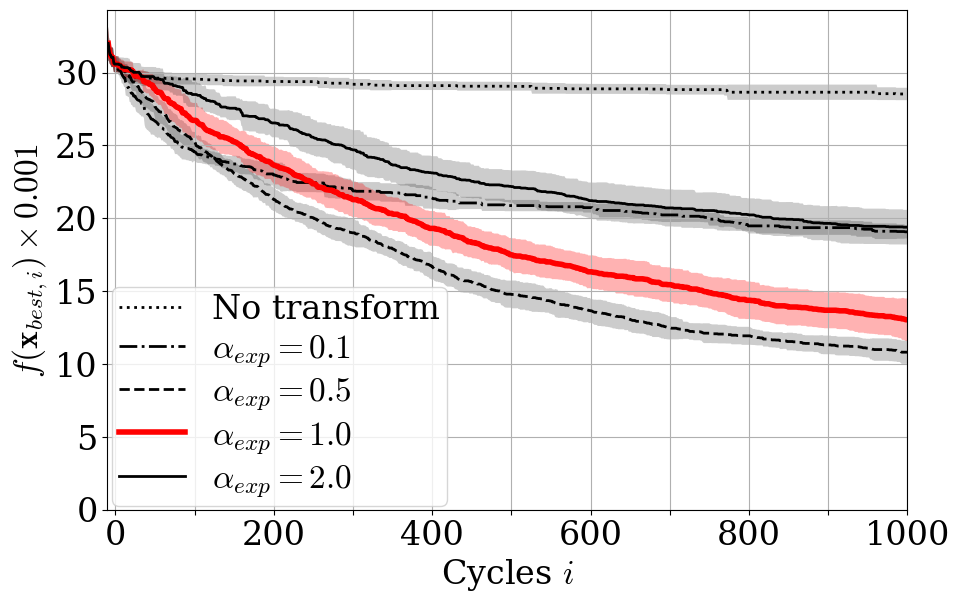}\label{fig:rosenbrock_binary_kq_d640_exp}}
    \mylab{-3.5cm}{-0.5cm}{(c)}
    \hspace{-0.4cm} 
    \subfloat{  \includegraphics[trim=-0.5cm 0cm 0cm 0cm, clip=true, height=4.4cm]{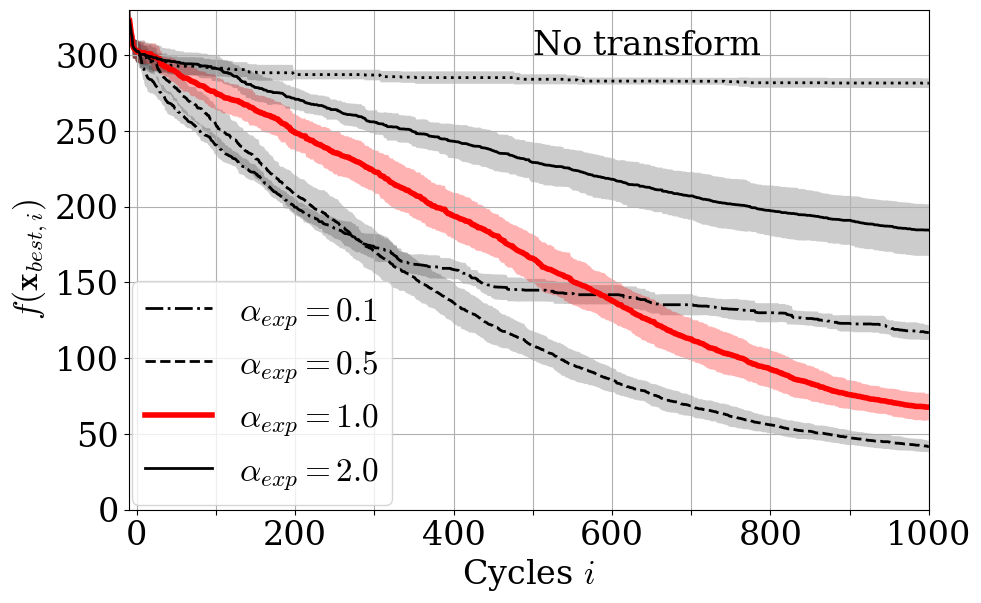}\label{fig:rastrigin_binary_kq_d640_exp}}
    \mylab{-3.5cm}{-0.5cm}{(d)} 
    \vspace{0.6cm}
  }
  \caption{Evolution of $f(\mathbf{x}_{best, i})$ with different $\alpha_{exp}$. (a, b) real variables with $d=80$ and (c, d) binary variables with $d=640$ for (a, c) Rosenbrock and (b, d) Rastrigin functions.}
\label{fig:transform}
\end{figure}

Sec.~\ref{sec:transform} described the exponential transformation of the training data ($f(\mathbf{x})$ output values) to help the surrogate model be constructed appropriately for optimization even when the black-box function yields an extensive dynamic range. The exponential transformation requires a transformation parameter $c_m$, which can be determined based on the initial training data, typically, the averaged output value $\alpha_{exp} \langle \mathbf{y}_{init}\rangle$ of the samples in the initial training dataset with $\alpha_{exp}$ being unity by default. This section discusses the effect of this transformation and the sensitivity of $c_m$ (or $\alpha_{exp}$) on the overall optimization performance.

Figure~\ref{fig:transform} shows the optimization history with kernel-QA under the conditions \texttt{r80} and \texttt{b640} for the Rosenbrock and Rastrigin functions. Four $\alpha_{exp}=0.1, 0.5, 1.0, 2.0$ and one condition without the exponential transformation are considered for each case. The results with $\alpha_{exp}=1.0$ are identical to the ones for \texttt{r80} and \texttt{b640} in Figs.~\ref{fig:d5_d80} and \ref{fig:d40_d640}. There are substantial improvements for all the cases using this transformation method with basically any $\alpha_{exp}$. The unity $\alpha_{exp}$ (default in the present study) seems quite a decent choice. However, the results show that twice or half of the value also yields reasonable (or better) optimization performance.

%%%%%%%%%%%%%%%%%%%%%%%%%%%%%%%%%%%%%%%
\subsection{Effect of initial data size}
\label{sec:init}

Figure~\ref{fig:typical_flow} shows that serial optimization methods require initial training data to construct the first model function. A relatively sizeable initial data size may result in a ``better'' model function for the first few optimization cycles. In contrast, a smaller initial data size reduces the cost of evaluating black-box functions during the initial data construction. Also, with a smaller initial data size, the ratio of newly added data samples during the optimization cycles to the total data becomes more prominent. This may result in faster convergence as, at the same cycles, the surrogate model fits better to the samples obtained during the optimization cycles rather than the initial data, for example, obtained randomly. The number of samples in the initial training data in the present study is set to be $b_{init}=10$ as summarized in Table~\ref{tab:conditions}. Here, the effect of initial data size on the optimization progress is discussed. 

\begin{figure}[!t]
\hspace{-0.5cm}
  \vspace{0.5cm}
  \centerline{
    \subfloat{  \includegraphics[trim=0cm 0cm 0cm 0cm, clip=true, height=4.4cm]{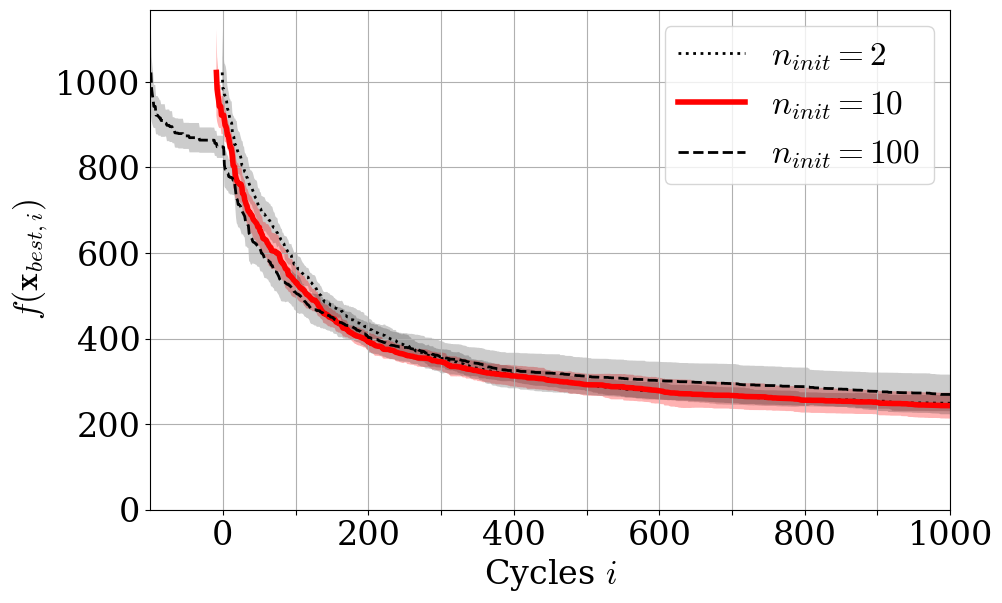}\label{fig:rastrigin_kq_d80_ninit}}
    \mylab{-3.5cm}{-0.5cm}{(a)} 
    \hspace{-0.4cm}
    \subfloat{  \includegraphics[trim=-0.5cm 0cm 0cm 0cm, clip=true, height=4.4cm]{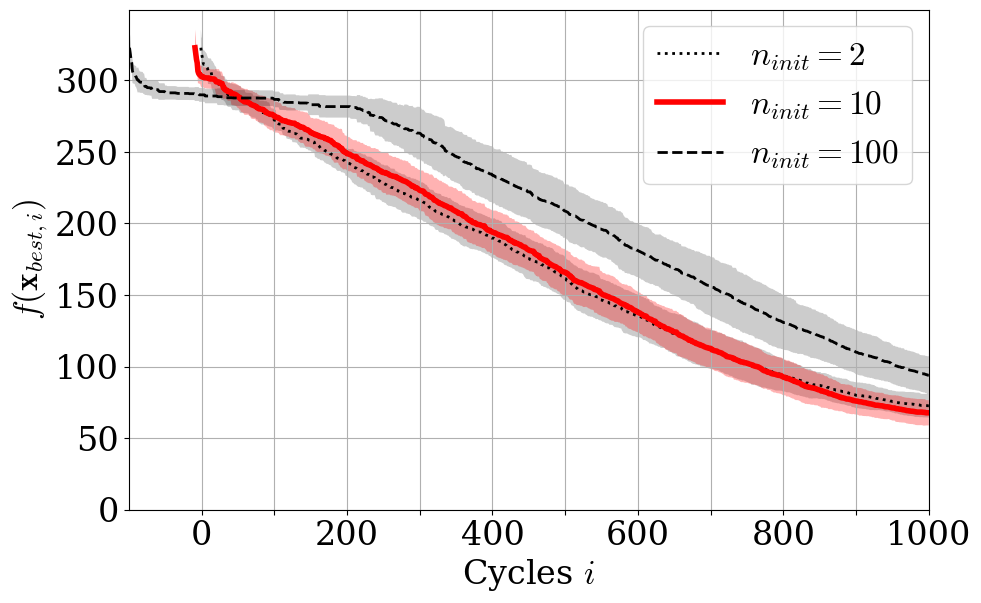}\label{fig:rastrigin_binary_kq_d640_ninit}}
    \mylab{-3.5cm}{-0.5cm}{(b)} 
  }
  \caption{Evolution of $f(\mathbf{x}_{best, i})$ with different initial data sizes $n_{init}=2$ (dotted), $10$ (red thick solid) and $100$ (dashed) for the Rastrigin function with kernel-QA: (a) \texttt{r80} and (b) \texttt{b640}.}
\label{fig:init}
\end{figure}

Figure~\ref{fig:init} shows optimization history with different initial data sizes $n_{init}=2$, $10$ and $100$. While Fig.~\ref{fig:init} shows the optimization results for the Rastrigin function, the trend for the corresponding Rosenbrock function is very similar. As clearly shown, the choice of $n_{init}$ does not substantially affect the overall optimization performance. That said, the result with $n_{init}=100$ shows lagged evolution of $f(\mathbf{x}_{best,i})$ reduction for the binary-variable problem in Fig.~\ref{fig:rastrigin_binary_kq_d640_ninit}. Also, such $n_{init}$ requires more black-box evaluations, which is impractical. While $n_{init} = 10$ is considered the default in the present experiments, the optimization cycles in kernel-QA can start with as small as $n_{init=2}$ samples without affecting optimization performance.

%
%
%

%%%%%%%%%%%%%%%%%%%%%%%%%%%%%%%%%%%%%%%
\subsection{Effect of $\beta$}
\label{sec:beta}

The assessments for kernel-QA hereinbefore utilized a surrogate model (Eq.~\eqref{eq:mean}), rather than an acquisition function (Eq.~\eqref{eq:acquisition}). However, for some BBO problems, it may be beneficial to consider both Eqs.~\eqref{eq:mean} and \eqref{eq:sigma} in the context of the acquisition function. The balance of the expected value (exploitation) and uncertainty (exploration) is controlled by $\beta$ in Eq.~\eqref{eq:acquisition}. when $\beta=0$, Eq.~\eqref{eq:acquisition} is the surrogate model.

\begin{figure}[!t]
  \vspace{0.5cm}
  \centerline{
    \hspace{-0.7cm}
    \subfloat{  \includegraphics[trim=0cm 0cm 0cm 0cm, clip=true, height=4.5cm]{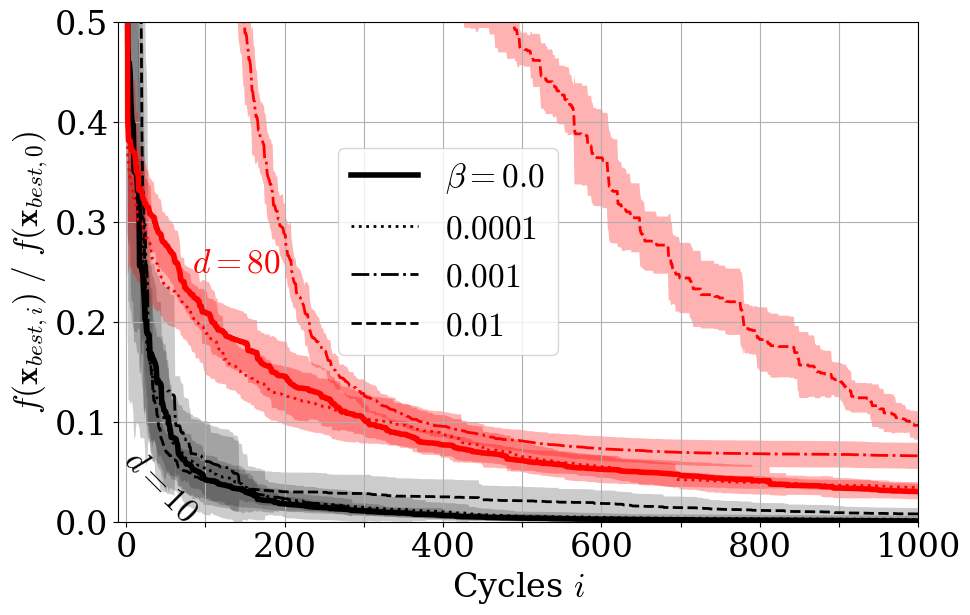}\label{fig:rosenbrock_kq_d80_beta}}
    \mylab{-3.5cm}{-0.5cm}{(a)}
    \hspace{-0.4cm}
    \subfloat{  \includegraphics[trim=0cm 0cm 0cm 0cm, clip=true, height=4.5cm]{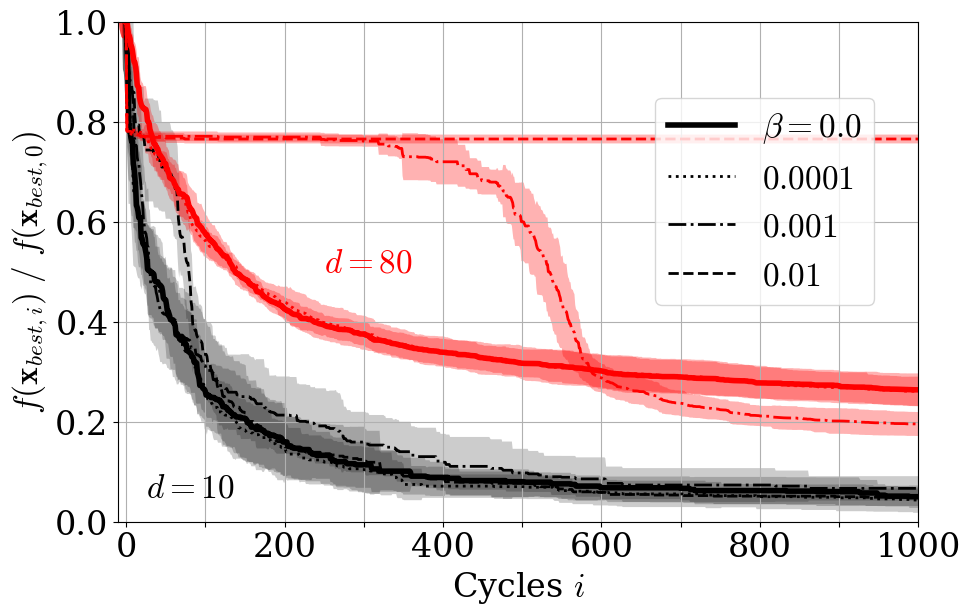}\label{fig:rastrigin_kq_d80_beta}}
    \mylab{-3.5cm}{-0.5cm}{(b)} 
  }
  \vspace{0.3cm}
  \centerline{
    \hspace{-0.6cm}
    \subfloat{  \includegraphics[trim=0cm 0cm 0cm 0cm, clip=true, height=4.5cm]{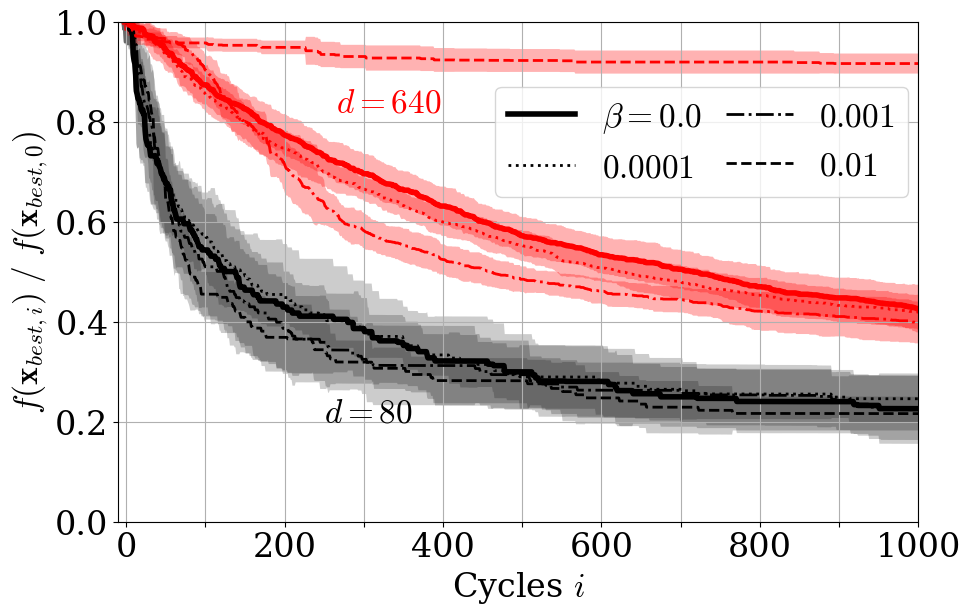}\label{fig:rosenbrock_binary_kq_d640_beta}}
    \mylab{-3.5cm}{-0.5cm}{(c)}
    \hspace{-0.4cm}
    \subfloat{  \includegraphics[trim=-0.5cm 0cm 0cm 0cm, clip=true, height=4.5cm]{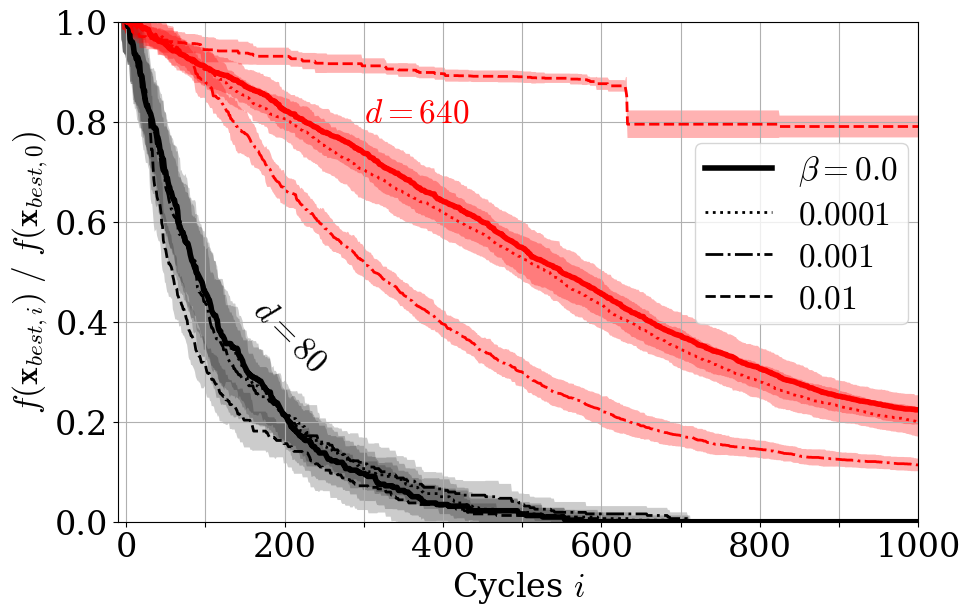}\label{fig:rastrigin_binary_kq_d640_beta}}
    \mylab{-3.5cm}{-0.5cm}{(d)} 
    \vspace{0.6cm}
  }
  \caption{Evolution of $f(\mathbf{x}_{best, i})$ with different ``exploration'' intensities (see Eq.~\eqref{eq:acquisition}), $\beta=0$ (thick), $0.0001$ (dotted), $0.001$ (dash-dotted) and $0.01$ (dotted). (a, b) \texttt{r10} and \texttt{r80}, and (c, d) \texttt{b80} and \texttt{b640} for (a, c) Rosenbrock and (b, d) Rastrigin functions. Optimization is performed by using kernel-QA. The shown value is normalized by using the best value at the beginning of the optimization cycles $f(\mathbf{x}_{best, 0})$.}
\label{fig:beta}
\end{figure}

Figure~\ref{fig:beta} compares evolution of $f(\mathbf{x}_{best, i})$ for the Rosenbrock and Rastrigin functions with zero and non-zero $\beta$ for the conditions \texttt{r10}, \texttt{r80}, \texttt{b80} and \texttt{b640} summarized in Table~\ref{tab:conditions}. The considered non-zero $\beta$ values are 0.0001, 0.001, and 0.01, where $\beta=0.01$ means more emphasis on exploration than $\beta=0$ (default value). The effect of $\beta$ seems more prominent, in both positive and negative manner, for the larger dimension, \texttt{r80} than \texttt{r10} for the real-variable case, and \texttt{b640} than \texttt{b80} for the binary-variable case.

For the more significant dimension cases (\texttt{r80} and \texttt{b640}), the optimization performance with $\beta=0.01$ (and larger $\beta$, expectedly) seems poorer regardless of functions or variable types. This result suggests the adverse effect of (excess) exploration for larger dimension problems, where the number of samples needed to sufficiently ``cover'' the search space grows rapidly with dimensionality. As for smaller non-zero $\beta$, its influence on overall optimization performance differs depending on the functions and variable types, and the trend does not seem consistent entirely. Such inconsistency is due to approximating the standard deviation in Eq.~\eqref{eq:sigma}. However, this inconsistency reduces towards the end of optimization cycles, and at the 1000-th cycle, the values of $f(\mathbf{x})$ converge to similar values between $\beta=0, 0.0001$, and $0.001$ for each case.

Although the present experiment does not demonstrate the consistent effectiveness of $\beta$, the result shows that the optimization performance could be further improved at optimization cycles $i\ll1000$ for relatively large dimension problems if an appropriate $\beta$ value can be used. Alternatively, the QUBO-compatible formulation of $\sigma$ (simplification from Eq.~\eqref{eq:sigma} to Eq.~\eqref{eq:raw_sigma}) could be explored to achieve consistent performance in future studies.

%
%
%

%%%%%%%%%%%%%%%%%%%%%%%%%%%%%%%%%%%%%%%
%%%%%%%%%%%%%%%%%%%%%%%%%%%%%%%%%%%%%%%
%%%%%%%%%%%%%%%%%%%%%%%%%%%%%%%%%%%%%%%
\section{Conclusions}
\label{sec:conclusions}

In the present study, a BBO method, kernel-QA, has been proposed, which is based on relatively low-order polynomial kernels and quadratic-optimization annealing, similar to FMQA. Instead of using a well-known acquisition function, kernel-QA considers a surrogate model function constructed analytically in a QUBO-compatible form to leverage the optimization using the Ising machine. Therefore, kernel-QA considers exploitation alone, whereas typical Bayesian optimization does both exploitation and exploration. Using a relatively low-order model function avoids model overfitting given relatively small training samples typical in BBO, and this feature helps circumvent local optimization.

Kernel-QA, Bayesian optimization, and partly FMQA have been assessed using artificial landscapes, the Rosenbrock and Rastrigin functions at various input dimensions (up to 80 for real and 640 for binary variables), variable types, and optimization conditions. For all test functions and variable types considered, kernel-QA performs well in terms of the final solution, the evolution of best objective function values, and its run-to-run standard deviation. The performance difference between the two methods becomes more transparent for larger-dimension problems or functions with local minima. The second point is especially explicitly clarified by comparing the history of the best objective value and distance between the best solution and true solution, and kernel-QA (as well as FMQA) showed a tendency to avoid lingering in local minima.

Several optimization parameters, exponential transformation coefficient ($\alpha_{exp}$) and size of initial data size ($n_{init}$), are also explored around their default values of $\alpha_{exp}=1$ and $n_{init}=10$. While $\alpha_{exp}=1$ seems reasonable, changing this parameter twice or half does not unduly influence the optimization performance. As for the initial dataset size, $n_{init}=2$, also showed almost identical optimization history. In contrast, we do not recommend using $n_{init}=100$ due to its negative influence on the performance and cost of constructing such a dataset.

Finally, a slight extension of kernel-QA to involve exploration is also proposed and assessed. An LCB parameter $\beta$ controls the balance of exploitation and exploration, and different $\beta$ values are considered here. The assessment shows the negative effect of exploration for a relatively large $\beta$ for the larger-dimension problems. Such an adverse effect is because the difficulty of sufficiently covering the search space grows rapidly, and the probability of sampling "promising" regions diminishes exponentially. Although the assessment showed some positive influence of exploration with appropriate $\beta$ values, the effect seems inconsistent, and this extension requires further improvement.

These results suggest that the proposed kernel-QA is a suitable BBO method for black-box functions with local minima and relatively large input dimensions.

\begin{appendices}

%%%%%%%%%%%%%%%%%%%%%%%%%%%%%%%%%%%%%%%
%%%%%%%%%%%%%%%%%%%%%%%%%%%%%%%%%%%%%%%
%%%%%%%%%%%%%%%%%%%%%%%%%%%%%%%%%%%%%%%
\section{Computational cost}
\label{sec:time}

\begin{figure}[!t]
  \centerline{
    \subfloat{  \includegraphics[trim=0cm 0cm 0cm 0cm, clip=true, height=5cm]{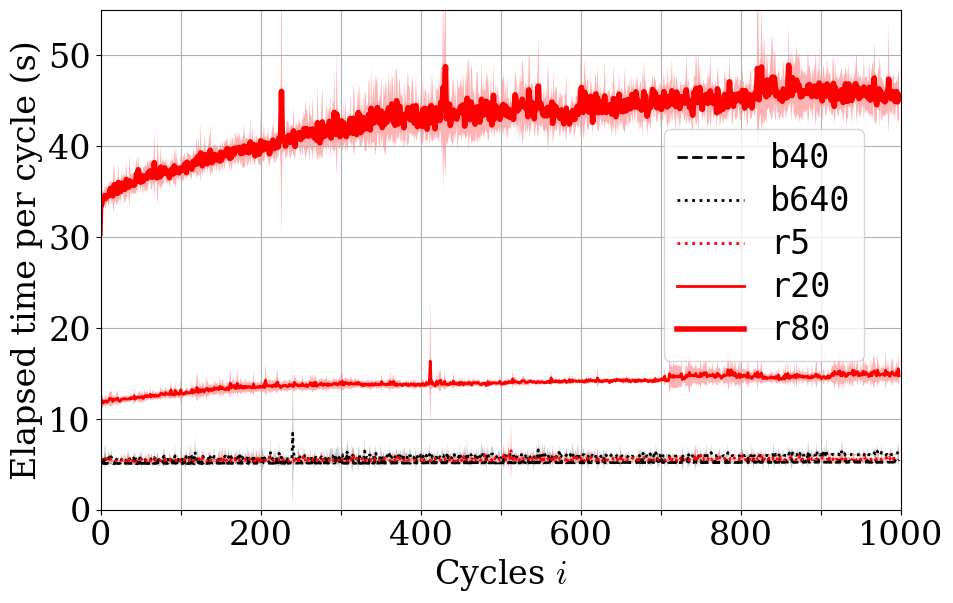}}
  }
  \caption{Variation of per-cycle computational time for kernel-QA for the binary-variable dimensions $d=40$ (black dashed) and $640$ (black dotted), and the real-variable dimensions $d=5$ (red dotted), $20$ (red thin solid) and $80$ (red thick solid) for the Rastrigin function. The plot shows the average and standard deviation values of ten independent runs.}
\label{fig:rastrigin_kq_time}
\end{figure}

The effect of problem dimensions and number of cycles on the computational time is discussed here. Figure~\ref{fig:rastrigin_kq_time} shows the variation of per-cycle computational time for kernel-QA under the conditions \texttt{b40}, \texttt{b640}, \texttt{r5}, \texttt{r20} and \texttt{r80} for the Rastrigin function. The per-cycle computational time is almost constant, except for \texttt{r80}, but the cycle dependency on computational time for \texttt{r80} does not seem strong (less than linear dependency).

As mentioned in Sec.~\ref{sec:conditions}, the annealing timeout for the optimization step ((2) in Fig.~\ref{fig:typical_flow}) is five seconds for all the cases. Thus, for the cases \texttt{b40}, \texttt{b640}, \texttt{r5}, annealing cost is predominant in overall computational cost. The computational time increases with the increase of the problem dimension (to be specific, the number of binary variables $d_B$ after variable conversion). Compared to \texttt{r5} with $d_B=300$ converted binary variables,  \texttt{r20} and \texttt{r80} shows an average of 2.5 and 7.8 times more computational costs, whereas their $d_B$ are 4 and 16 times. Most of the cost increase is due to the QUBO formulation and preprocessing of request data before annealing on the Ising machine.

%
%
%

%%%%%%%%%%%%%%%%%%%%%%%%%%%%%%%%%%%%%%%
%%%%%%%%%%%%%%%%%%%%%%%%%%%%%%%%%%%%%%%
%%%%%%%%%%%%%%%%%%%%%%%%%%%%%%%%%%%%%%%
\section{Exponential transformation for Bayesian optimization}

\begin{figure}[!t]
  \centerline{
    \subfloat{  \includegraphics[trim=0cm 0cm 0cm 0cm, clip=true, height=5cm]{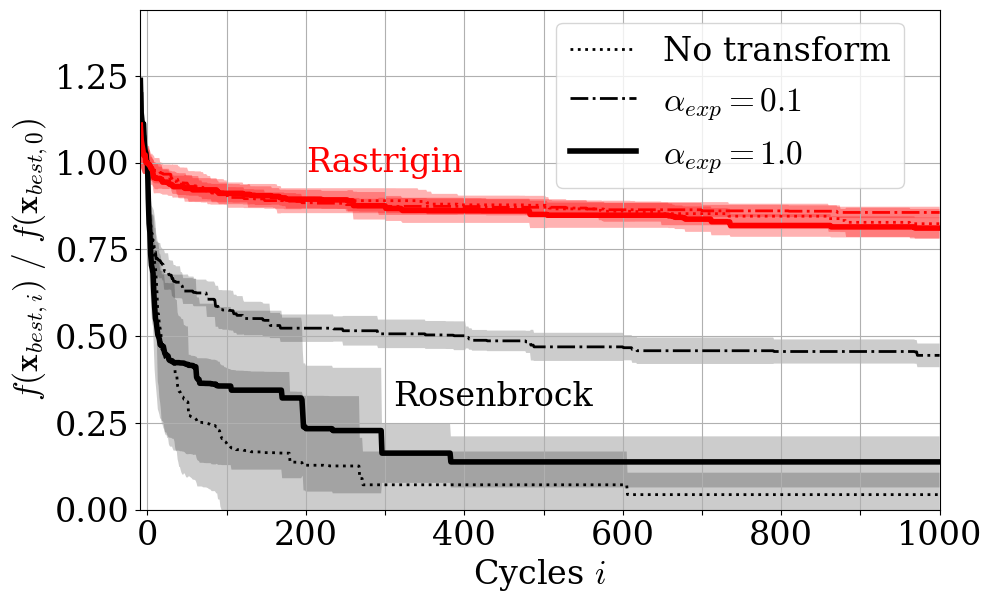}}
  }
  \caption{Evolution of $f(\mathbf{x}_{best, i})$ with different conditions for the exponential transformation (Sec.~\ref{sec:transform}) for Bayesian optimization under the condition \texttt{r80}. The shown value is normalized by using the best value at the beginning of the optimization cycles $f(\mathbf{x}_{best, 0})$.}
\label{fig:bo_d80_exp}
\end{figure}

The exponential transformation described in Sec.~\ref{sec:transform} and assessed in Sec.~\ref{sec:a_exp} is a robust method to facilitate the surrogate model construction in kernel-QA. This transformation is applied to Bayesian optimization, and the results are compared in Fig.~\ref{fig:bo_d80_exp}.

The effect of the exponential transformation on the Bayesian optimization results is not as substantial as kernel-QA shown in Sec.~\ref{sec:a_exp}. Also, the method performs best without the transformation. The result implies the predominance of exploration in this optimization method in large-dimension problems.

\end{appendices}

\newpage
\clearpage

%%%%%%%%%%%%%%%%%%%%%%%%%%%%%%%%%%%%%%%
%%%%%%%%%%%%%%%%%%%%%%%%%%%%%%%%%%%%%%%
%%%%%%%%%%%%%%%%%%%%%%%%%%%%%%%%%%%%%%%
\section*{Declarations}

\subsection*{Author contributions}

All authors contributed to the study conception and design. Material preparation, data collection and analysis were performed by Yuki Minamoto. The draft of the manuscript was written by all authors. All authors read and approved the final manuscript.

\subsection*{Data availability}

The optimization program code used in the present study is available from \url{https://pypi.org/project/amplify-bbopt/0.1.0/#files}. The resulting data from the code for the conditions presented in the study is accessible upon request from Yuki Minamoto (yuki.minamoto@fixstars.com)

\subsection*{Funding}

No funding was received to assist with the preparation of this manuscript.

\subsection*{Conflict of interest}

The authors declare that they have no conflict of interest.

\bibliography{sn-bibliography}% common bib file
%% if required, the content of .bbl file can be included here once bbl is generated
%%\input sn-article.bbl

\end{document}